\documentclass[11pt,a4paper]{article}
\usepackage{amssymb,amsmath,amsthm,amsfonts}
\usepackage[top=2.5cm, bottom=3cm]{geometry}

\usepackage{authblk}

\usepackage{hyperref}
\hypersetup{
    colorlinks=true,       
    linkcolor=black,          
    citecolor=green,        
    filecolor=magenta,      
    urlcolor=blue           
}

\usepackage{stmaryrd}
\usepackage{subfig}
\usepackage{tikz}

\usetikzlibrary{arrows}

\usepackage{float}

\def\N{\mathbb{N}}

\def\R{\mathbb{R}}
\def\C{\mathbb{C}}

\def\I{\mathbb{I}}

\def\D{{\cal D}}
\def\LL{{\cal L}}

\def\CPP{C\kern-.05em\raise.23ex\hbox{\footnotesize{+}\kern-.07em\footnotesize{+}}}

\def\d{\delta}
\def\e{\varepsilon}

\def\h{\lambda}
\def\a{\alpha}

\def\dim{\mathrm{dim}}

\newcommand\norm[1]{\interleave#1\interleave}
\newcommand\x[1]{\mathsf{#1}}
\newcommand\ex[2][]{\langle #2 \rangle_{#1}}
\newcommand\co[2][]{\,\rangle\hspace{-1pt}#2\hspace{-1pt}\langle_{#1}}

\newcommand{\interval}[1]{\left\llbracket #1 \right\rrbracket}

\theoremstyle{plain}
\newtheorem{theorem}{Theorem}[section]
\newtheorem{lemma}[theorem]{Lemma}
\newtheorem{proposition}[theorem]{Proposition}
\newtheorem{observation}[theorem]{Observation}

\newtheorem{remark}[theorem]{Remark}
\newtheorem{corollary}[theorem]{Corollary}

\theoremstyle{definition}
\newtheorem{definition}[theorem]{Definition}
\newtheorem{example}[theorem]{Example}

\def\Item$#1${\item $\displaystyle#1$   \hfill\refstepcounter{equation}{\normalfont(\theequation)}}

\definecolor{jszary}{gray}{0.5}

\title{\bf On rigorous estimates of eigenspaces and eigenvalues of a matrix}

\author{\L{}ukasz Struski\thanks{Corresponding author}}
\author{Jacek Tabor}
\author{Piotr Zgliczy\'nski \footnote{Research has been supported by Polish
National Science Centre grant 2011/03B/ST1/04780}  }

\affil{
\small Jagiellonian University,\\
\small Faculty of Mathematics and Computer Science,\\
\small \L{}ojasiewicza 6, 30-348 Krak\'ow, Poland\\\vspace{5pt}
\small e-mail: \href{mailto:struski@im.uj.edu.pl}{\sf struski@im.uj.edu.pl},\\ \hspace{25pt}\href{mailto:tabor@ii.uj.edu.pl}{\sf tabor@ii.uj.edu.pl}, \\ \hspace{45pt}\href{mailto:umzglicz@cyf-kr.edu.pl}{\sf umzglicz@cyf-kr.edu.pl}
}

\date{\today}

\begin{document}
 \maketitle

\begin{abstract}
We present a method of cones for rigorous estimations of eigenvectors, eigenspaces and eigenvalues of a matrix.
The key notion is the cone-domination and is inspired by ideas from hyperbolic dynamical systems. We present theorems
which allow to rigorously locate the spectrum of the matrix and the eigenspaces, also multidimensional ones in case
of eigenvalues of multiplicity greater than one or clusters of close eigenvalues. 

In case of isolated eigenvalue we show that the our method give the same or better estimates than ones known in literature.

\end{abstract}

\bigskip

\noindent\small {\bf Key words and phrases:} eigenvectors, eigenvalues, Gerschgorin theorem, cone condition, spectrum of the matrix\\[0.2em]
\small {\bf 2010 Mathematics Subject Classification:} 65F15, 37D30.

\medskip

\section{Introduction}

Eigenvalues and eigenvectors are the basic tools used in mathematics and computer science (linear algebra, differential equations, statistics, etc.). Currently, there are a lot of numerical methods, which allow to solve (not rigorously) the eigenproblem \cite{golubmatrix, watkinsfundamentals, watkinsmatrix}. 
However, in such fields as computer assisted proofs for PDEs, methods that allow us to specify the rigorous bounds on the eigenvalues (see \cite{wilczak}) are required with increasing interest. 
In papers \cite{saad, Stewart_2} we can find perturbation theory for eigenvalues and invariant subspaces of matrices, but presented there techniques are not very useful for these problems. 
In this paper we present tools to find rigorous bounds for eigenvalues (all or some of them) and their corresponding eigenspaces.

Assume that we have a square matrix $A$ which the entries (or blocks) on the diagonal 'dominate' the off-diagonal entries (blocks) and we want to obtain efficient computable bounds (a formula) for the spectrum and eigenspaces of $A$.
Regarding the bounds on the spectrum almost all of the known methods  are given by the Gerschgorin theorems and its modifications, for example the Brauer ovals \cite{Br,Va} or the generalization of the Gerschgorin theorem to the case of  multi-dimensional blocks by Feingold and Varga \cite{Gen}.
Estimation of isolated eigenvectors from Gerschgorin's results are due to Stewart \cite{Stewart} and Wilkinson \cite{Wilk}. However, Wilkinson's result does not give the whole eigenspace
in the case of not simple eigenvalue or a cluster of close eigenvalues.
In  \cite{Ya}, T. Yamamoto showed how find rigorous error bounds for computed single eigenvalues and eigenvectors of real matrices on the basis of an existence theorem for solutions of nonlinear systems using iteration Newton's method.
However,  the Yamamoto's approach gives no theoretical estimates for the bounds for computed eigenvalues and eigenvectors.

In this article we propose  a new method for the estimates of eigenvalues and eigenspaces.
Our approach is based on the ideas coming from the hyperbolic dynamics \cite{Nh}
and can be illustrated by the following simple two-dimensional example.

\begin{example}
 Consider the matrix $A$ is defined by the formula
\(
A=
\begin{bmatrix}
0 & 2 \\
1 & 4
\end{bmatrix}.
\)
Note that if we take the gray cone (see Figure \ref{Ob:1}) and we start to iterate points of this cone by matrix $A$ then range of our cone will be reduced to the eigenspace corresponding to one of the eigenvectors of $A$.

\begin{figure}[H]
  \centering
\begin{tikzpicture}[
    scale=0.25,
    axis/.style={thin, ->, >=stealth'},
    important line/.style={thick},
    dashed line/.style={dashed, thin},
    every node/.style={color=black}
    ]
\begin{scope}[]
\path [fill=gray] (-5,5) -- (0,0) -- (5,5);
\path [fill=gray] (-5,-5) -- (0,0) -- (5,-5);
\path (0,-6) node[below] {\scriptsize $\{x=(x_1,x_2)\ :\ |x_1|\leq|x_2|\}$};
\draw[axis] (-5,0)  -- (5.3,0) node(xline)[right]{};
\draw[axis] (0,-5) -- (0,5.2) node(yline)[above] {};
\node [right](ref1) at (5,5){};
\end{scope}

\begin{scope}[xshift=15cm]
\path [fill=gray] (-3.34,-5) -- (0,0) -- (-2,-5);
\path [fill=gray] (3.34,5) -- (0,0) -- (2,5);
\draw[axis] (-5,0)  -- (5.3,0) node(xline)[right]{};
\draw[axis] (0,-5) -- (0,5.2) node(yline)[above] {};
\node [left](ref2) at (-4,5){};
\node [left](ref3) at (6,5){};
\end{scope}

\begin{scope}[xshift=30cm]
\path [fill=black] (-2.3,-5) -- (0,0) -- (-2.2,-5);
\path [fill=black] (2.2,5) -- (0,0) -- (2.3,5);
\draw[thin, gray] (-2.247,-5) -- (0,0) -- (2.247, 5);
\draw[axis] (-5,0)  -- (5.3,0) node(xline)[right]{};
\draw[axis] (0,-5) -- (0,5.2) node(yline)[above] {};
\node [left](ref4) at (-4,5){};
\node [left](ref5) at (6,5){};
\end{scope}

\begin{scope}[xshift=45cm]
\draw[thin, gray] (-2.247,-5) -- (0,0) -- (2.247, 5);
\draw[axis] (-5,0)  -- (5.3,0) node(xline)[right]{};
\draw[axis] (0,-5) -- (0,5.2) node(yline)[above] {};
\node [left](ref6) at (-4,5){};
\end{scope}

\draw[red, dashed, very thick,->] (ref1) .. node[above] {$\pmb{A}$} controls (7,6) and (9,6) ..  (ref2);
\draw[red, dashed, very thick,->] (ref3) .. node[above] {$\pmb{A}$} controls (22,6) and (24,6) ..  (ref4);
\draw[red, dashed, very thick,->] (ref5) .. node[above] {$\pmb{A}$} controls (37,6) and (39,6) ..  (ref6);
\end{tikzpicture}
  \caption{Transformation of the cone by the matrix $A$. \label{Ob:1}}
\end{figure}
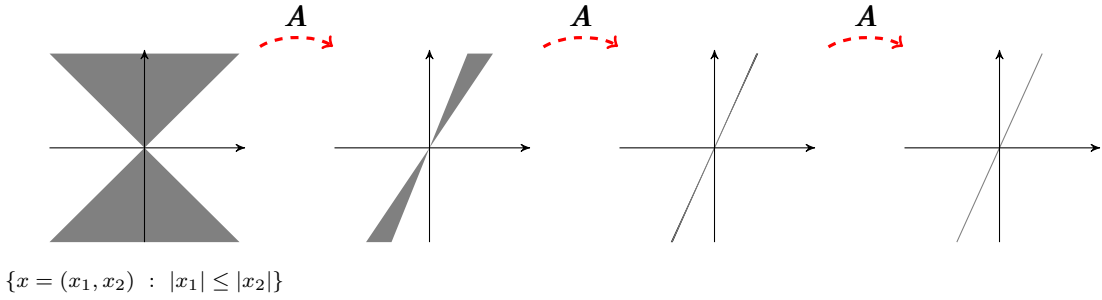

\noindent Iterating backward ($A^{-1}$) the cone $\{x=(x_1,x_2)\ :\ |x_1|\geq|x_2|\}$ we obtain second eigenspace.
\end{example}

This example illustrates that we can estimate eigenspace by using an invariant  cone. We started to study this problem and it turned out that using forward and backward invariant cones we were able to give not only good bounds for eigenvectors but also for eigenvalues. In addition, by means of this tool we can locate the eigenspaces and eigenvalues of products of many matrices.

To explain our main results we introduce some basic notations.
Let
$\|\x{x}\|:=\max\limits_j |x_j|$. For  $\x{x}=(x_1,\ldots, x_k,\ldots,
x_n)\in\R^n$ we set
\[
\|\x{x}\|_{\leq k}=\max\limits_{i\leq
k}|x_i|  \quad\mbox{ and }\quad   \|\x{x}\|_{>k}=\max\limits_{i> k}|x_i|.
\]
For linear map
$A\colon\R^k\times\R^{n-k}\to\R^k\times\R^{n-k}$ we define
the extension and contraction  constants:
\begin{align*}
\co{A} & = \;\inf \{ R \in \R_+ \, | \, \|A\x{x}\|\leqslant R\cdot\|x\| \text{ for all } \x{x}\in\R^n: \|A\x{x}\|_{\leq k}\geq\|A\x{x}\|_{>k}\}, \\[0.3em] 
\ex{A} & =  \sup \{ R \in \R_+ \, | \, \|A\x{x}\| \geqslant
R\cdot\|x\|  \text{ for all } \x{x}\in\R^n: \|\x{x}\|_{\leq k}\leq
\|\x{x}\|_{>k}\}.
\end{align*}
Observe that these constants can be obtained  by the optimal solution of a standard constrained optimization problem.

We say that $A$ is {\em dominating} if $\co{A}\, <\ex{A}$.
It turns out that composition of dominating maps is dominating, see Proposition~\ref{proposition:1}.
Now we are ready to present two main results in our paper

\medskip

\noindent {\bf Main Result I [Theorem \ref{thm:gerssh-dominating}].}{\em \/
Let $A\in\C^{n\times n}$ be a matrix with an isolated Gerschgorin disk. Then $A$ is dominating.
}

\medskip

\noindent Together with the following result we get that our method is generalization of  the Gerschgorin theorem in the case of the isolated Gerschgorin disk of multiplicity one.

\medskip

\noindent {\bf Main Result II [simplified version of Theorem
\ref{thm:cone-main-eigen-location}].}{\em \/ Let
$A\colon\R^k\times\R^{n-k}\to\R^k\times\R^{n-k}$ be dominating.
Then there exists a unique direct sum decomposition $F_1\oplus
F_2=\R^n$ into $A$-invariant subspaces $F_1$, $F_2$ such that
\[
\sigma(A|_{F_1})\subset\overline{B}(0,\co{A}),\quad
\sigma(A|_{F_2})\subset\C\setminus B(0,\ex{A}).
\]

Moreover, we have:
\begin{enumerate}
\item $\dim F_1=k$,\; $\dim F_2=n-k$,
\item $F_1\subset\{\x{x}\in\R^n : \|\x{x}\|_{\leq k}\geq \|\x{x}\|_{>k}\}$ \quad and \quad
$F_2\subset\{\x{x}\in\R^n : \|\x{x}\|_{\leq k}\leq
\|\x{x}\|_{>k}\}$,
\item $\|A|_{F_1}\|\leq \co{A}\quad \mbox{and}\quad \|(A|_{F_2})^{-1}\|\leq \ex{A}^{-1}$.
\end{enumerate}
}

\noindent In comparison with  the Gerschgorin's theorems our method has the following advantages:
\begin{itemize}
\item locate spectrum and eigenspaces of a matrix when multiple eigenvalues or clusters of very close eigenvalues are present,
\item gives better estimation for isolated eigenvalues,
\item allow to deal with composition of matrices.
\end{itemize}

\smallskip

The content of this paper can be briefly described as follows: in Section \ref{sec:cones} we introduce notion of cones and build the concept of dominating matrix.  In Section \ref{sec:main-res} we establish the main result: Theorem \ref{thm:cone-main-eigen-location} which allow us to rigorously estimate eigenspaces and eigenvalues.
In Section \ref{sec:eigen-estm} we develop computable estimates for the eigenvalues and eigenspaces based on the results from the Section \ref{sec:main-res}. In Section~\ref{sec:comparison} we compare the proposed method with the Gerschgorin theorem in the case of the isolated Gerschgorin disk. We show that all matrices which have an isolated Gerschgorin disk, are dominating and if the radius of this disk nonzero, we  obtain sharper bounds.  This means that our approach can be used whenever  the Gerschgorin disk is isolated.
We also show examples of matrices for which we can not use the Gerschgorin theorem since the Gerschgorin disks cannot be separated, but our method still works, see Example~\ref{ex:mA-better-G}.

Now we present some basic notation.
By $\R$ and $\C$ we denote the sets of real, and
complex numbers.
The spectrum $\sigma(A)$ of a square matrix
$A=[a_{ij}]\in\C^{n\times n}$ we define the collection of all
eigenvalues of $A$, i.e.
\[
\sigma(A):=\{\h\in\C\; :\; A-\h I \text{ is singular}\}.
\]

By $I_{\C^n}$ we mean the identity matrix of size $n$, while $\I$ denotes the interval
$\interval{-1,1}$. For $\e>0$ we put
$B_{\C}(0,\e):=\{z\in\C\;:\; |z|<\e\}$.

\section{Cones and dominating maps}
\label{sec:cones}

 In this section we introduce the basic concepts
and tools  of our  method of invariant cones to locate the eigenspaces and bound the spectrum for  matrices.
For this end we modify the concept of cones from \cite{KT}.
 Our approach is strongly motivated  by the methods from the theory of hyperbolic dynamical systems, in particular by  the results of Newhouse \cite{Nh}, who obtained conditions for hyperbolic splitting on compact invariant set for a diffeomorphism in terms of its induced action on a cone-field and its complement.

\begin{definition}
By a {\em cone-space} we understand a finite dimensional Banach
space $E$ with semi-norms $\co[]{\cdot}$ (we call it {\em
contracting}), $\ex{\cdot}$ (which we call {\em expanding}) such
that
\[
\norm{\x{x}} := \max(\co[]{\x{x}}, \ex{\x{x}})
 \]
 defines an equivalent norm on $E$.
 By the {\em r-norm} for $r>0$ on  the cone-space $E$ we take
\[
   \norm{\x{x}}_r:=\max(\co[]{\x{x}}, r\cdot\ex{\x{x}}).
\]
\end{definition}

\begin{definition}\label{def:1}
Let $E$ be a cone-space. We define the {\em $r$-contracting cone}
in $E$ by
\begin{align*}
   \co[r]{E} & :=\{\x{x} \in E: \;\co[]{\x{x}} \geq r\ex{\x{x}}\},\\
   \intertext{and the {\em $r$-expanding cone} in $E$ by}
   \ex[r]{E} & :=\{\x{x} \in E:\; \co[]{\x{x}} \leq r\ex{\x{x}} \}.
\end{align*}
\end{definition}
\noindent Note that
\begin{equation}\label{eq:1}
E= \co[r]{E}\cup\ex[r]{E}.
\end{equation}
In the same way we define $r$-contracting cone and $r$-expanding
cone in subspace of $E$. If $r=1$ we will omit the subscript $r$,
in particular we speak of contracting cone. We introduce the
scaling by $r$ of semi-norms to have a better control over size of
the cones (see Figure \ref{rys:1}), which will consequently allow
us to better locate the eigenvectors.
\begin{figure}[H]
  \centering
    \subfloat[The contracting cone in $\R\times\R$.]{
\begin{tikzpicture}[
    scale=0.45,
    axis/.style={thin, dashed, ->, >=stealth'}
    ]
\path [fill=jszary] (-5,5) -- (0,0) -- (-5,-5);
\path [fill=jszary] (5,5) -- (0,0) -- (5,-5);
\draw[axis] (-5,0)  -- (5.4,0);
\draw[axis] (0,-5) -- (0,5.4);
\end{tikzpicture}
}
        \hspace{2cm}
\subfloat[The $2$-contracting cone in $\R\times\R$.]{
\begin{tikzpicture}[
    scale=0.45,
    axis/.style={thin, dashed, ->, >=stealth'}
    ]
\path [fill=jszary] (-5,2.5) -- (0,0) -- (-5,-2.5);
\path [fill=jszary] (5,2.5) -- (0,0) -- (5,-2.5);
\draw[axis] (-5,0)  -- (5.4,0);
\draw[axis] (0,-5) -- (0,5.4);
\end{tikzpicture}
}
\end{figure}

\begin{figure}[H]
\ContinuedFloat
\setcounter{figure}{2}
  \centering
\subfloat[ The expanding cone in $\R\times\R$.]{
\begin{tikzpicture}[
    scale=0.45,
    axis/.style={thin, dashed, ->, >=stealth'}
    ]
\path [fill=jszary] (-5,5) -- (0,0) -- (5,5);
\path [fill=jszary] (-5,-5) -- (0,0) -- (5,-5);
\draw[axis] (-5,0)  -- (5.4,0);
\draw[axis] (0,-5) -- (0,5.4);
\end{tikzpicture}
         }
        \hspace{2cm}
\subfloat[The $2$-expanding cone in $\R\times\R$.]{
\begin{tikzpicture}[
    scale=0.45,
    axis/.style={thin, dashed, ->, >=stealth'}
    ]
\path [fill=jszary] (-5,5) -- (-5,2.5) -- (0,0) -- (5,2.5) -- (5,5) -- (-5,5);
\path [fill=jszary] (-5,-5) -- (-5,-2.5) -- (0,0) -- (5,-2.5) -- (5,-5) -- (-5,-5);
\draw[axis] (-5,0)  -- (5.4,0);
\draw[axis] (0,-5) -- (0,5.4);
\end{tikzpicture}
}
\caption{The cones in the cone-space $\R\times\R$.}\label{rys:1}
\end{figure}
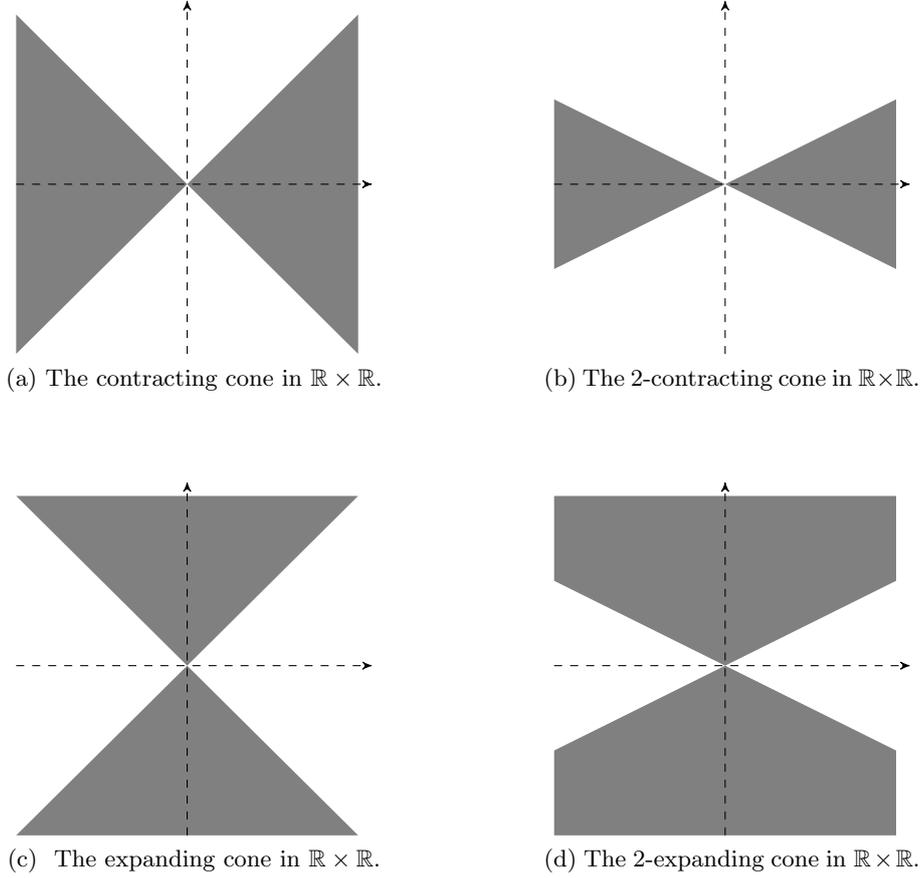
If $E$ has a fixed product structure $E = E_1\times E_2$, we
introduce a natural cone-space structure on $E$ by defining
seminorms
\[
 \co[]{\x{x}}:=\|x_1\|, \quad \ex{\x{x}}:=\|x_2\|\quad\text{ for }\;\x{x}=(x_1,x_2)\in E_1\times E_2.
\]

\medskip

In the proof of our main result, Theorem \ref{thm:cone-main-eigen-location}, the following proposition
will play a crucial role.

\begin{proposition}\label{prop:1}
Let $E=E_1\times E_2$ be a cone-space and let $r>0$ be given.
Assume that we have direct sum decomposition $E=V_1\oplus V_2$
such that
\[
V_1\subset\co[r]{E}\quad\text{and}\quad V_2\subset\ex[r]{E}.
\]
Then \( \dim V_1=\dim E_1\quad\text{and}\quad\dim V_2=\dim E_2. \)
\end{proposition}
\begin{proof}
Let $n:=\dim E_1$ and $m:=\dim E_2$. First we show that $\dim
V_1\leq n$. For an indirect proof, assume that $\dim V_1>n$. Then
there exist linearly independent vectors
$\x{v_1},\ldots,\x{v_{n+1}}\in V_1$. Obviously $\x{v_i}=(w_i,
z_i)$ for $i\in\{1,\ldots,n+1\}$ and unique $w_i\in E_1$, $z_i\in
E_2$. Since $w_1,\ldots,w_{n+1}\in E_1$ and $\dim E_1=n$ there
exist a set of $n+1$ scalars, $\a_1, \ldots,\a_{n+1}$, not all
zero, such that
\[
\a_1 w_1+\ldots+\a_{n+1}w_{n+1}=0.
\]
Note that
\[
\x{z}:=\a_1 z_1+\ldots+\a_{n+1}z_{n+1}\neq 0,
\]
because otherwise the vectors $\x{v_1},\ldots,\x{v_{n+1}}$ would
not be linearly independent. Consequently we obtain
\begin{equation*}
(0,\x{z})=\left(\sum_{i=1}^{n+1}\a_i w_i, \sum_{i=1}^{n+1}\a_i
z_i\right)\in V_1\subset\co[r]{E},
\end{equation*}
and thus $r\|\x{z}\|\leq\|0\|$, which implies that $\x{z}=0$. We
get a contradiction with the fact the sequence of vectors
$\x{v_1},\ldots,\x{v_{n+1}}$ is linearly independent.

The proof that $\dim V_2\leq m$ is analogous. Finally, since $\dim
E=n+m$ and $\dim V_1\leq n$, $\dim V_2\leq m$ we obtain
\[
\dim V_1=n,\quad\text{and}\quad\dim V_2=m.
\]
\end{proof}

By an {\em operator} we mean a linear mapping between cone-spaces
$E$ and $F$. We denote the space of all operators by $\LL(E,F)$.
If $F=E$, we denote $\LL(E,E)$ by $\LL(E)$.

\medskip

Let $A\in\LL(E,F)$. We define
\begin{alignat}{2}
   \co[r]{A} & := & \inf & \{R \in \R_+ \, | \,  \norm{A\x{x}}_r \leqslant R\norm{\x{x}}_r \text{ for all } \x{x}\in E: A\x{x} \in \co[r]{F}\}, \label{eq:2}\\[0.3em]
   \ex[r]{A} & :=\; &\sup &\{ R \in \R_+ \, | \,  \norm{A\x{x}}_r \geqslant R\norm{\x{x}}_r \text{ for all }  \x{x} \in E: \x{x} \in \ex[r]{E}\}. \label{eq:3}
\end{alignat}

The following lemma is obvious.
\begin{lemma}
\label{lem:coA-exA-alter}
Let $A\in\LL(E,F)$.
\begin{align}
\co[r]{A} &=  \left(\inf \{ \norm{\x{x}}_r \, | \, A\x{x} \in \co[r]{F}, \, \norm{ A\x{x}}_r=1 \}\right)^{-1} \; \mbox{when $A$ is invertible},  \label{eq:coA-alter} \\[0.3em]
   \ex[r]{A} &= \inf \{ \norm{A\x{x}}_r \, | \, \x{x} \in \ex[r]{E}, \, \norm{\x{x}}_r=1  \}.  \label{eq:exA-alter}
\end{align}
\end{lemma}

\begin{remark}\label{rem:rate}
Observe, that
\begin{alignat*}{2}
\norm{A\x{x}}_r &\leqslant \co[r]{A}\norm{\x{x}}_r &\mbox{for}\quad & \x{x}\in A^{-1}\co[r]{F},\\[0.3em]
\norm{A\x{x}}_r &\geqslant \ex[r]{A}\norm{\x{x}}_r &\quad\mbox{for}\quad & \x{x}\in \ex[r]{E}.
\end{alignat*}
\end{remark}

\noindent The above definitions of $
\co[r]{A}$ and  $\ex[r]{A}$ are modifications of analogous notions
in \cite{Nh}, where  $\ex[]{A}$ is called  the expansion rate and
$1/\!\!\co[]{A}$ is  the co-expansion rate. Using of those rates
we can generalize the classical dominating maps which are relevant
to our research.

\begin{definition}\label{def:2}
We say that $A\in\LL(E,F)$ is {\em $r$-dominating}, if
   \[
   \co[r]{A} < \ex[r]{A}.
   \]
By $\D_r(E,F)$ we denote the set of all $A\in\LL(E,F)$ which are
$r$-dominating. If $F=E$, we denote the space $\D_r(E,E)$ by
$\D_r(E)$.
\end{definition}

\begin{observation}\label{ob:1}
Let $\tilde{E}\subset E$, $\tilde{F}\subset F$ be subspaces and
let $A\in\LL(E,F)$ be such that $A(\tilde{E})\subset\tilde{F}$.
Then $A|_{\tilde E} \in \LL(\tilde E,\tilde F)$ and
\[
\co[r]{A|_{\tilde{E}}}\leq\co[r]{A}\quad\text{ and
}\quad\ex[r]{A}\leq\ex[r]{A|_{\tilde{E}}}.
\]
Moreover, if $A\in\D_r(E,F)$ then $A\in\D_r(\tilde{E},\tilde{F})$.
\end{observation}
\begin{proof}
It is a consequence of \eqref{eq:2}, \eqref{eq:3} and Definition
\ref{def:2}.
\end{proof}

It turns out that $r$-cones are invariant for $r$-dominant operators.

\begin{theorem}\label{tw:1}
Let $A\in \D_r(E,F)$ and let  $\x{v}\in E$ be arbitrary. Then
\begin{align*}
  \x{v}\in \ex[r]{E} & \implies A\x{v}\in \ex[r]{F} ,\\[0.3em]
  A\x{v}\in \co[r]{F} & \implies \x{v}\in \co[r]{E}.
\end{align*}
\end{theorem}

\begin{proof}
The proof is a simple modification of the proof of
\cite[Proposition 2.1]{KT}.
\end{proof}

As a consequence of the above theorem we obtain that composition of $r$-dominating maps is
$r$-dominating. Moreover, we get estimate for expansion and contraction rates.

\begin{proposition}\label{proposition:1}
Let $A\in \D_r(F,G)$ and $B\in \D_r(E,F)$. Then $A\circ B\in
\D_r(E,G)$ and
\begin{equation}\label{eq:4}
  \co[r]{A\circ B}\leq \co[r]{A}\cdot\co[r]{B},\ \ex[r]{A\circ B} \geq \ex[r]{A}\cdot\ex[r]{B}.
\end{equation}
\end{proposition}

\begin{proof}
To prove the first inequality from \eqref{eq:4}, consider an
$\x{x}\in E$ such that $(A\circ B)(\x{x})\in \co[r]{G}$. From
\eqref{eq:2} and Theorem \ref{tw:1} we know that $B\x{x}\in
\co[r]{F}$, and thus we have
\[
  \norm{A\circ B(\x{x})}_r\leq\co[r]{A} \cdot \norm{B\x{x}}_r
\leq \co[r]{A} \cdot \co[r]{B} \cdot \norm{\x{x}}_r.
\]
Hence
\[
  \co[r]{A\circ B} \leq \co[r]{A} \cdot \co[r]{B}.
\]
Using \eqref{eq:3} and Theorem \ref{tw:1}, we obtain the second
inequality from \eqref{eq:4}.

As a simple consequence of \eqref{eq:4} we obtain $A\circ B\in
\D_r(E,G)$.
\end{proof}

In the remainder of this section we show
how to estimate $\co[r]{A}$, $\ex[r]{A}$.
Consider two cone-spaces $E = E_1\times E_2$  and $F = F_1\times
F_2$. Let $A\colon E\to F$ be an operator given in the matrix
form by
\[
 A = \begin{bmatrix}
       A_{11} &A_{12}           \\[0.3em]
       A_{21} & A_{22}
     \end{bmatrix}.
\]
 By
\[
\norm{A} _{r}:= \max\big(\|A_{11}\| + \frac{1}{r}\|A_{12}\|,
r\|A_{21}\| + \|A_{22}\|\big)
\]
we define the {\em $r$-norm} of operator $A$, where $\|.\|$ is an
operator norm. Observe that it satisfies
\[
\norm{A\x{x}} _{r}\leq \norm{A} _{r}\cdot \norm{\x{x}} _r \quad
\text{ for }\; \x{x}\in E.
\]
Note that in general it is not (except for the case when $E_1$ is
one dimensional) the  operator  norm for $\norm{\cdot}_{r}$.

\begin{theorem} \label{thm:formulas-contr-exp}
Let $A=[A_{ij}]_{1\leq i,j\leq 2}\in\LL(E_1\times E_2,F_1\times
F_2)$ and $r\in(0,\infty)$ be given.
\begin{enumerate}
\item We have
\[
\co[r]{A} \leq \| A_{11}\|+\frac{1}{r}\|A_{12}\|.
\]
\item Additionally, if $A_{22}$ is invertible, then
\[
\ex[r]{A} \geq \|A_{22}^{-1}\|^{-1} - r\|A_{21}\|.
\]
\end{enumerate}
\end{theorem}
\begin{proof}
For the proof of the first inequality, we take $\x{x}=(x_1,x_2)\in
E_1\times E_2$ such that $A\x{x}\in \co[r]{F}$. From Definition
\ref{def:1} we have
\begin{equation}\label{eq:6}
\|A_{11}x_1+A_{12}x_2\|\geq r\|A_{21}x_1+A_{22}x_2\|,
\end{equation}
and therefore
\begin{align*}
\norm{A\x{x}}_r &= \max(\|A_{11}x_1+A_{12}x_2\|, r\|A_{21}x_1+A_{22}x_2\|)\\
&\overset{by\,\eqref{eq:6}}{=}\|A_{11}x_1+A_{12}x_2\|\leq \|A_{11}\|\cdot \|x_1\| + \frac{1}{r}\|A_{12}\|\cdot r\|x_2\|\\
&\leq (\|A_{11}\| + \frac{1}{r}\|A_{12}\|)\cdot \norm{\x{x}}_r.
 \end{align*}

For the proof of the second inequality, suppose that
$\x{x}=(x_1,x_2)\in \ex[r]{E}$, where $x_1\in E_1$, $x_2\in E_2$.
Then
\begin{equation}\label{eq:7}
\|x_1\|\leq r\|x_2\|=\norm{\x{x}}_r.
\end{equation}
We know that
\begin{equation}\label{eq:8}
\|A_{22}x_2\|\geq \|A_{22}^{-1}\|^{-1}\|x_2\| \geq 0.
\end{equation}
Finally, we obtain
\begin{align*}
\norm{A\x{x}}_r &\geq r\|A_{21}x_1+A_{22}x_2\|\geq r\|A_{22}x_2\| - r\|A_{21}x_1\|\\
&\overset{by\,\eqref{eq:8}}{\geq}r\|A_{22}^{-1}\|^{-1}\|x_2\| -
r\|A_{21}\|\|x_1\|\overset{by\,\eqref{eq:7}}{\geq}\left(\|A_{22}^{-1}\|^{-1}
-  r\|A_{21}\|\right)\cdot\norm{\x{x}}_r.
\end{align*}

\end{proof}

\begin{example}
\label{ex:dominat-ineq}
Let us verify that the matrix $A\in\LL(\C\times\C,\C\times\C)$, \(
 A =
 \begin{bmatrix}
  2 & 1.5 \\
  1 & 5
 \end{bmatrix}
\) is dominating. By Theorem \ref{thm:formulas-contr-exp} we have \( \co[]{A} \leq
3.5  < 4 \leq\ex{A}, \) and therefore $A$ is dominating.
\end{example}

Let us stress that the estimates from Theorem~\ref{thm:formulas-contr-exp} are sharp, but there are cases when we do not have equalities in them.

\begin{example}
\label{ex:dominat-exactly}
Let $A\in\LL(\C\times\C,\C\times\C)$ be given by the formula \(
 A =
 \begin{bmatrix}
  2 & 3 \\
  2 & 5
 \end{bmatrix}.
\) We show that $A$ is dominating. Observe that  Theorem \ref{thm:formulas-contr-exp} does not allow us to decide    whether this matrix $A$ is dominating, since
\[
\co[]{A}\leq \| A_{11}\|+\|A_{12}\|=5 \ \mbox{ and }\ \ex[]{A}\geq \|A_{22}^{-1}\|^{-1} - \|A_{21}\|=3.
\]
We calculate exactly $\co{A}$ and $\ex{A}$ (we take the norm $\|\cdot\|_{\infty}$) from the formulas \eqref{eq:coA-alter} and \eqref{eq:exA-alter}. The minimum of \eqref{eq:exA-alter} is realized in points $(1,-1)^T$ and $(-1,1)^T$. It is easy to see that the matrix $A$ is invertible, so minimum of \eqref{eq:coA-alter} is realized in points $(\frac{1}{2},0)^T$ and $(-\frac{1}{2},0)^T$.
Hence
\[
\co{A}=2 \quad\mbox{ and }\quad \ex{A}=3.
\]
Finally, we obtain that $A$ is dominating.
Observe that for this example Gerschgorin theorem does not hold (it is impossible to separate Gerschgorin disks).
\end{example}

\section{Localization of eigenspaces based on cones and dominating maps}
\label{sec:main-res} In this section we show that the eigenspaces of the
$r$-dominating operator $A$ lie in the corresponding $r$-cones.
Moreover, we can estimate $\sigma(A)$ with
the help of $\co[r]{A}$, $\ex[r]{A}$.

\begin{lemma}\label{lem:1}
Let $A\in\D_r(E)$. Then
\begin{equation}\label{eq:9}
\h\in\sigma(A) \implies |\h|\in
[0,\!\!\co[r]{A}]\cup[\ex[r]{A},\infty).
\end{equation}
Moreover $[0,\!\!\co[r]{A}]\cap [\ex[r]{A},\infty)=\emptyset$.
\end{lemma}

\begin{proof}
Since $A\in\D_r(E)$ we get $[0,\!\!\co[r]{A}]\cap
[\ex[r]{A},\infty)=\emptyset$.

Now we show implication \eqref{eq:9}. Let $\h$ be an eigenvalue of
$A$ and let $\x{x}\in E$ be a corresponding eigenvector. By
\eqref{eq:1} we know that $\x{x}\in\co[r]{E}\cup\ex[r]{E}$. We
consider two cases. First suppose that $\x{x}\in\co[r]{E}$. Since
$\x{x}$ is an eigenvector, $A\x{x} = \h\x{x}$, and thus
$A\x{x}\in\co[r]{E}$. By \eqref{eq:2} we get
\[
|\h|\leq\co[r]{A}.
\]
Now suppose that $\x{x}\in\ex[r]{E}$. By \eqref{eq:3} we get
\[
|\h|\geq\ex[r]{A},
\]
which completes the proof.
\end{proof}

Let $E$ be a finite dimensional vector space over the field $\C$
and let operator $A\colon E\to E$ be given.  One can easily deduce from the Jordan theorem (see also \cite[Appendix to Chapter
4]{irwin} for the general case) that if $\sigma(A) =
\sigma_1\cup\sigma_2$  then there is a unique direct sum
decomposition $E=E_{\sigma_1}\oplus E_{\sigma_2}$ such that
$A(E_{\sigma_1})\subset E_{\sigma_1}$, $A(E_{\sigma_2})\subset
E_{\sigma_2}$ and $\sigma(A|_{E_{\sigma_1}})=\sigma_1$,
$\sigma(A|_{E_{\sigma_2}})=\sigma_2$. For $0<c<d$ we define
\[
E_{\leq c}:=E_{\{\h\; :\; |\h|\leq c\}}\quad\text{and}\quad
E_{\geq d}:=E_{\{\h\; :\; |\h|\geq d\}}.
\]

\begin{theorem}\label{thm:r-dom-spectrum-gap}
Let $E$ be a finite dimensional cone-space and let $A\in\D_r(E)$.
Then there is a direct sum decomposition
$E=E_{\leq\co[r]{A}}\oplus E_{\geq\ex[r]{A}}$ which satisfies
\[
E_{\leq\co[r]{A}} \subset  \co[r]{E},\; E_{\geq\ex[r]{A}} \subset
\ex[r]{E}.
\]
\end{theorem}
\begin{proof}
From Lemma \ref{lem:1} and the comments preceding our theorem we
obtain a decomposition of $E$ into $A$-invariant subspaces
\[
E=E_{\leq\co[r]{A}}\oplus E_{\geq\ex[r]{A}},
\]
such that
\[
\sigma(A|_{E_{\leq\co[r]{A}}})=\{\lambda :
|\lambda|\in[0,\!\!\co[r]{A}]\}\quad\text{and}\quad
\sigma(A|_{E_{\geq\ex[r]{A}}})=\{\lambda :
|\lambda|\in[\ex[r]{A},\infty)\}.
\]

Now we show $E_{\leq\co[r]{A}} \subset  \co[r]{E}$. Consider an
arbitrary $\x{x}\in E_{\leq\co[r]{A}}$. The case when $\x{x}=0$ is
obvious. Assume that $\x{x}\neq 0$. Without any loss of the generality we can assume that $\|\x{x}\| = 1$. For an
indirect proof, assume that $\x{x}\notin\co[r]{E}$. Then
by (\ref{eq:1}) we get $\x{x}\in\ex[r]{E}$. Let $\e>0$ be arbitrary. From
the fact that $\x{x}\in E_{\leq\co[r]{A}}$, we know that
\begin{equation}\label{radius}
\limsup\limits_{m\rightarrow +\infty}\sqrt[m]{\norm{A|_{E_{\leq\co[r]{A}}}^m}}=\sup\sigma(A|_{E_{\leq\co[r]{A}}})\leq\co[r]{A}.
\end{equation}
Note that inequality \eqref{radius} holds for all norms.
For all $x\in E_{\leq\co[r]{A}}$ we obtain
\[
\limsup\limits_{m\rightarrow +\infty}\sqrt[m]{\norm{A^m
\x{x}}}\leq\co[r]{A},
\]
and thus there exists an $M\in\N$ such that for all $m\in\N$
\[
m\geq M \Rightarrow \sqrt[m]{\norm{A^m \x{x}}}\leq \co[r]{A}+\e.
\]
Since $\x{x}\in\ex[r]{E}$ and from Theorem \ref{tw:1} we obtain
\[
\x{x}\in\ex[r]{E} \Rightarrow A\x{x}\in \ex[r]{E} \Rightarrow  \cdots
\Rightarrow A^m\x{x}\in\ex[r]{E}.
\]
Using \eqref{eq:3} and Remark \ref{rem:rate} we get
\begin{align*}
\norm{A\x{x}}&\geq \ex[r]{A} \norm{\x{x}}, \\
\norm{A^2\x{x}} = \norm{A(A\x{x})}&\geq \ex[r]{A} \norm{A\x{x}}\geq \ex[r]{A}^2 \norm{\x{x}}, \\
&\;\, \vdots\\
\norm{A^m\x{x}}&\geq \ex[r]{A}^m \norm{\x{x}}.
\end{align*}
Finally we have
\[
\ex[r]{A} = \sqrt[m]{\ex[r]{A}^m} \leq \sqrt[m]{\norm{A^m \x{x}}} \leq
\co[r]{A}+\e.
\]
Since $\e$ was arbitrary, we get a contradiction with the fact
that $A$ is $r$-dominating.

Analogously, to prove inclusion $E_{\geq\ex[r]{A}}\subset\ex[r]{E}$,
assume that $\x{x}\in E_{\geq\ex[r]{A}}$ and $\x{x}\notin\ex[r]{E}$.
Then  $\x{x}\in\co[r]{E}$. Since
$\sigma(A|_{E_{\geq\ex[r]{A}}})=\sigma_{{\geq\ex[r]{A}}}:=\{\h\; :\;
|\h|\geq\ex[r]{A}\}$ and $0\notin\sigma_{{\geq\ex[r]{A}}}$ we know that
$A|_{E_{\geq\ex[r]{A}}}\colon E_{\geq\ex[r]{A}}\to E_{\geq\ex[r]{A}}$ is
invertible. Let $\e>0$ be arbitrary. Using the fact that $\x{x}\in
E_{\geq\ex[r]{A}}$, by dual result \eqref{radius}, we know that
\[
 \limsup\limits_{m\rightarrow +\infty}\sqrt[m]{\norm{A|_{ E_{\geq\ex[r]{A}}}^{-m} \x{x}}}\leq\ex[r]{A}^{-1},
 \]
and thus there exists an $M\in\N$ such that for all $m\in\N$
\begin{equation}\label{eq:10}
m\geq M \Rightarrow\sqrt[m]{\norm{A|_{ E_{\geq\ex[r]{A}}}^{-m}
\x{x}}}\leq \ex[r]{A}^{-1}+\e.
\end{equation}
From the Observation \ref{ob:1} and Theorem \ref{tw:1} we get
\[
\x{x}\in\co[r]{ E_{\geq\ex[r]{A}}} \Rightarrow A|_{
E_{\geq\ex[r]{A}}}^{-1}\x{x}\in \co[r]{ E_{\geq\ex[r]{A}}} \Rightarrow
\cdots  \Rightarrow A|_{ E_{\geq\ex[r]{A}}}^{-m}\x{x}\in\co[r]{
E_{\geq\ex[r]{A}}},
\]
and from \eqref{eq:2} and Remark \ref{rem:rate} we have
\begin{align*}
\norm{\x{x}}&\leq \co[r]{A|_{ E_{\geq\ex[r]{A}}}} \norm{A|_{ E_{\geq\ex[r]{A}}}^{-1}\x{x}}, \\
\norm{A|_{ E_{\geq\ex[r]{A}}}^{-1}\x{x}} &\leq \co[r]{A|_{ E_{\geq\ex[r]{A}}}} \norm{A|_{ E_{\geq\ex[r]{A}}}^{-2}\x{x}}, \\
&\;\, \vdots\\
\norm{A|_{ E_{\geq\ex[r]{A}}}^{-m+1}\x{x}}&\leq \co[r]{A|_{
E_{\geq\ex[r]{A}}}} \norm{A|_{ E_{\geq\ex[r]{A}}}^{-m}\x{x}}.
\end{align*}
Hence
\begin{equation}\label{eq:11}
\norm{\x{x}}\leq (\co[r]{A|_{ E_{\geq\ex[r]{A}}}})^m\norm{A|_{
E_{\geq\ex[r]{A}}}^{-m}\x{x}}.
\end{equation}

Finally from the Observation \ref{ob:1} and \eqref{eq:10},
\eqref{eq:11} we obtain
\[
\co[r]{A}\geq\co[r]{A|_{ E_{\geq\ex[r]{A}}}} = \sqrt[m]{(\co[r]{A|_{
E_{\geq\ex[r]{A}}}})^{m}}\geq\sqrt[m]{\frac{1}{\norm{A|_{
E_{\geq\ex[r]{A}}}^{-m}
\x{x}}}}\geq\frac{1}{\ex[r]{A}^{-1}+\e}=\ex[r]{A}\cdot
\frac{1}{1+\e\cdot\ex[r]{A}},
\]
which gives a contradiction with the fact that $A$ is
$r$-dominating.
\end{proof}

Now we are ready to state the main result on the eigenspaces and eigenvalue location using
our method of cones and dominating maps.

\begin{theorem}  \label{thm:cone-main-eigen-location} 
Let $E=E_1\times E_2$ be a finite dimensional cone-space and let
$A\in\D_r(E)$. Then there exists a unique direct sum decomposition
$E=F_1\oplus F_2$ of $A$-invariant subspaces $F_1$, $F_2$ such
that
\[
\sigma(A|_{F_1})\subset\overline{B}(0,\co[r]{A}),\quad
\sigma(A|_{F_2})\subset\C\setminus B(0,\ex[r]{A}).
\]

Moreover, we have:
\begin{enumerate}
\item $\dim F_1=\dim E_1$,\; $\dim F_2=\dim E_2$,
\item $F_1\subset\co[r]{E}$\; and\; $F_2\subset\ex[r]{E}$,
\item $\|A|_{F_1}\|\leq \co[r]{A}\quad \mbox{and}\quad \|(A|_{F_2})^{-1}\|\leq \ex[r]{A}^{-1}$. \label{lab:1}
\end{enumerate}
\end{theorem}
\begin{proof}
From Theorem \ref{thm:r-dom-spectrum-gap} we know that exists  a unique direct sum
decomposition $E=E_{\leq\co[r]{A}}\oplus E_{\geq\ex[r]{A}}$ which
satisfies
\[
E_{\leq\co[r]{A}} \subset  \co[r]{E},\; E_{\geq\ex[r]{A}} \subset
\ex[r]{E}.
\]
We take $F_1=E_{\leq\co[r]{A}}$ and $F_2=E_{\geq\ex[r]{A}}$. By
Proposition \ref{prop:1} we obtain $\dim F_1=\dim E_1$ and  $\dim
F_2=\dim E_2$.

Now we show that
$\sigma(A|_{F_1})\subset\overline{B}(0,\co[r]{A})$. Let $\x{x}\in
F_1$ be an eigenvector of $A$ and let $\h$ be the eigenvalue of
$A$ corresponding to $\x{x}$. Since $\x{x}$ is an eigenvector
($A\x{x} = \h\x{x}$) and $F_1 \subset\co[r]{E}$ therefore
$A\x{x}\in\co[r]{E}$. By \eqref{eq:2} we obtain that
$|\h|\leq\co[r]{A}$, so we get
$\sigma(A|_{F_1})\subset\overline{B}(0,\co[r]{A})$.

Now suppose that $\x{x}\in F_2$. Since $F_2\subset\ex[r]{E}$ and
by \eqref{eq:3} we get $|\h|\geq\ex[r]{A}$. Hence $
\sigma(A|_{F_2})\subset\C\setminus B(0,\ex[r]{A})$.

The inequalities of item \ref{lab:1} we obtain from \eqref{eq:2}
and \eqref{eq:3}.
\end{proof}

As a direct consequence of the above theorem we obtain  the following conclusion.

\begin{corollary}\label{cor:2}
Let $r\in(0,\infty)$ and $n\in\N$. Assume that an operator
$A\in\D_r(\C\times\C^{n-1})$ is given. Then there
exists unique eigenvalue $\h$ of $A$ such that $|\h|\leq\co[r]{A}$
and  the eigenspace corresponding to $\h$ is one-dimensional. The
unique (after rescaling) eigenvector $\x{x}$ corresponding to the
eigenvalue $\h$ satisfies
\[
\x{x}\in(1,0,\ldots,0)^T+\{0\}\times
\overline{B}_{\C}(0,1/r)^{n-1}\subset(1,0,\ldots,0)^T+\frac{1}{r}\cdot(0,\I,\ldots,\I)^T
+ \frac{1}{r}\cdot(0,\I,\ldots,\I)^Ti.
\]
\end{corollary}
\begin{proof}
It is a direct consequence of Theorem  \ref{thm:cone-main-eigen-location} and Definition
\ref{def:1}.
\end{proof}

Because at the origin of our approach based on cones and dominating maps is the
theory of hyperbolic dynamical systems, so our method should be well suited to locate the eigenspaces and eigenvalues
of products of many matrices. In the example below we contrast our approach  with a naive approach, which tries
to diagonalize a matrix obtained as a product of many matrices. The essential feature of this example is that the matrices we multiply
are known with some accuracy only.

\begin{example}\label{example:iteration}
Let the matrices $A_i\in\LL(\R\times\R)$, $i\in\{1,\ldots,15\}$ be
such that
\[
  A_i\in
\begin{bmatrix}
\interval{0, 0.5} & \e\I \\
\e\I & \interval{1.5,2}
\end{bmatrix},
\]
where $\e=0.01$ and $\I=\interval{-1,1}$. Consider the matrix
$B:=A_{15}\cdot\ldots\cdot A_1$.

From Theorem \ref{thm:formulas-contr-exp} we obtain that $A_i\in\D(\R\times\R)$ and
\[
\co{A_{i}}\leq 0.5+\e, \ \ex{A_i}\geq 1.5-\e.
\]
From Theorem \ref{tw:1} and Proposition \ref{proposition:1} we
conclude that $B\in\D(\R\times\R)$ and
\[
\co{B}\leq \co{A_{15}}\cdot\ldots\cdot\co{A_1}, \quad
\ex{B}\geq\ex{A_{15}}\cdot\ldots\cdot\ex{A_1}.
\]
From Theorem \ref{thm:cone-main-eigen-location} we obtain that eigenvalues $\h_1$ and
$\h_2$ of $B$ such that
\[
|\h_1| \leq (0.5+\e)^{15} \ \text{ and }\ |\h_2| \geq
(1.5-\e)^{15}.
\]

Now, a naive method will ask first for a computation of $B$. Using interval arithmetic we obtained
\[
B\in \left[
\begin{array}{cc}
 \interval{-1.45687,1.45693} & \interval{-218.543,218.544} \\
 \interval{-218.543,218.544} & \interval{433.611,32782.94}
\end{array}
\right].
\]
However, there exists  matrix $B_1$ within the bounds given above,
which has both eigenvalues larger than $1$. For example, let us consider
\[
B_1= \left[
\begin{array}{cc}
 1 & 100 \\
 -100 & 521 \\
\end{array}
\right].
\]
This matrix have the eigenvalues $\h_1=21$ and $\h_2=501$.
Consequently, this means that none of the methods applied to the
product matrix will not give us the expected estimation
$|\lambda_1|<1$ and $|\lambda_2|>1$.
\end{example}

\section{Estimations of the eigenvalues and eigenvectors }
\label{sec:eigen-estm}

In this section we develop
computable estimates for the eigenvalues and eigenspaces based on
the results from the previous section.

\begin{lemma}\label{lem:2}
Let $A\in\LL(E_1\times E_2)$ be given such that
\[
A:=
\begin{bmatrix}
 A_{11} & A_{12}\\
 A_{21} & A_{22}
\end{bmatrix}.
\]
If $A_{22}$ is invertible, $d=\|A_{22}^{-1}\|^{-1}-\|A_{11}\|>0$
and $\Delta:=d^2-4\|A_{12}\|\|A_{21}\|>0$ then
\[
A\in\D_{r}(E_1\times E_2)\quad\text{for }\; \left\{
\begin{alignedat}{2}
&r\in\left(\frac{d-\sqrt{\Delta}}{2\|A_{21}\|}, \frac{d+\sqrt{\Delta}}{2\|A_{21}\|}\right) &&\text{ if }\;\|A_{21}\|\neq 0\\
&r\in\left(\frac{\|A_{12}\|}{d},\infty\right)&&\text{ if
}\;\|A_{21}\|=0
\end{alignedat}
\right..
\]
\end{lemma}
\begin{proof}
Let $a:=\|A_{12}\|$, $b:=\|A_{11}\|$ and $c:=\|A_{21}\|$. Making
use of Theorem \ref{thm:formulas-contr-exp} it suffices to show that
\[
b+\frac{a}{r}<(d+b)-cr.
\]
Multiplying both sides of the above inequality by the positive
number $r$ we get the inequality
\begin{equation}\label{eq:12}
cr^2-dr+a<0.
\end{equation}

If $c=0$ then we get $r>\frac{a}{d}$. Suppose now that  $c\neq
0$. Since from our assumption follows that $\Delta>0$ we see
inequality \eqref{eq:12} is satisfied for
\[
r\in\left(\frac{d-\sqrt{\Delta}}{2c},
\frac{d+\sqrt{\Delta}}{2c}\right).
\]
\end{proof}

\begin{remark}
\label{rem:good-eps} Let $A$ be an operator, which satisfies the
assumptions of Lemma \ref{lem:2} (in particular $\Delta>0$). Let  $a:=\|A_{12}\|$,
$b:=\|A_{11}\|$ and $c:=\|A_{21}\|\not = 0$.
It is easy to see, that
\[
    \frac{d-\sqrt{\Delta}}{2c} < \frac{d}{2c} < \frac{d+\sqrt{\Delta}}{2c} < \frac{d}{c}.
\]
\end{remark}

Therefore, if $A$ satisfies the
assumptions of Lemma \ref{lem:2} and $\|A_{21}\|\not = 0$ and we want to find possibly largest $r$ for which $A$ is $r$-dominating, then we can
take $r=\frac{d}{2\|A_{21}\|}$.  With this choice we have $r<r_{max}< 2r$, where $r_{max}$ is the supremum the set of $r$'s obtained in the above lemma,
therefore we might not be optimal, but we obtain easily manageable expression.


As a corollary we obtain a well-known result for the location of an isolated eigenvalue and its eigenspace \cite[Theorem 3.11]{Stewart}. We present its statement in the form adapted to our notation.

\begin{theorem} (\cite[Theorem 3.11]{Stewart}) \label{thm:cone-single-eigenval}
Let $A=[a_{ij}]_{1\leq i,j\leq n}\in\LL(\C\times\C^{n-1})$ be
given in the block from by
\[
A:=
\begin{bmatrix}
 A_{11} & A_{12}\\
 A_{21} & A_{22}
\end{bmatrix},
\]
where $A_{11}=a_{11}$. Assume that $A_{22}-a_{11}\cdot
I_{\C^{n-1}}$ is invertible and $\|(A_{22}-a_{11}\cdot
I_{\C^{n-1}})^{-1}\|^{-2}-4\|A_{12}\|\|A_{21}\|>0$. Then
\begin{enumerate}
\item there exists a unique eigenvalue $\h$ of $A$ which satisfies
 \[
|\h-a_{11}|\leq
2\|A_{12}\|\cdot\|A_{21}\|\cdot\|(A_{22}-a_{11}\cdot
I_{\C^{n-1}})^{-1}\|,
\]
\item the eigenspace corresponding to $\h$ is one-dimensional and there exist unique $\d_2$, $\ldots$, $\d_n\in\C$,
\[
\|(0,\d_2,\ldots,\d_n)^T\|\leq2\|A_{21}\|\cdot\|(A_{22}-a_{11}\cdot
I_{\C^{n-1}})^{-1}\|\cdot\|(1,0,\ldots,0)^T\|
\]
such that \( (1,\d_2,\ldots,\d_n)^T \) is the eigenvector
corresponding to $\h$.
\end{enumerate}
\end{theorem}
\begin{proof} It is easy to see that if $A_{21}=0$, then theorem holds. Therefore we will assume that $\|A_{21}\| >0$.

In order to apply Lemma~\ref{lem:2} to matrix $A - a_{11}I_{\mathbb{C}^n}$ we set
 $a:=\|A_{12}\|$, $c:=\|A_{21}\|$ and $d=\|(A_{22}-a_{11}\cdot
I_{\C^{n-1}})^{-1}\|^{-1}$.
By Lemma \ref{lem:2} and
Remark~\ref{rem:good-eps} we get $A-a_{11}\cdot
I_{\C^{n}}\in\D_{d/(2c)}(\C\times\C^{n-1})$, and from Corollary \ref{cor:2} and Theorem~\ref{thm:formulas-contr-exp} we conclude that
there exists a unique eigenvalue $\h$ of $A$ which satisfies
\begin{equation*}
  |\lambda - a_{11}| < \co[\frac{d}{2c}]{A-a_{11}I} \leq  \frac{1}{\frac{d}{2c}} \|A_{12}\| = \frac{2ac}{d}.
\end{equation*}

From Theorem \ref{thm:cone-main-eigen-location} (second point) we know  that eigenspace, which contains eigenvector corresponding to the $\lambda$, lies in $\co[\frac{d}{2c}]{\C\times\C^{n-1}}$. Hence (see Definition \ref{def:1}) we obtain unique $\d_2$, $\ldots$, $\d_n\in\C$, \(
\|(0,\d_2,\ldots,\d_n)^T\|\leq2\|A_{21}\|\cdot\|(A_{22}-a_{11}\cdot
I_{\C^{n-1}})^{-1}\|\cdot\|(1,0,\ldots,0)^T\| \) such that \(
(1,\d_2,\ldots,\d_n)^T \) is the eigenvector corresponding to
$\h$.
\end{proof}

Let us stress here that in the proof Theorem~\ref{thm:cone-single-eigenval} through Lemma~\ref{lem:2} we used estimates for $\co[]{A}$ and $\ex[]{A}$ provided by Theorem~\ref{thm:formulas-contr-exp}, which  may fail establish that a matrix is dominating for a dominating matrix. If this is
the case we will use Theorem~\ref{thm:cone-main-eigen-location}. This happens in Examples \ref{ex:G-better} and \ref{ex:not-G-better}.

The following lemma shows how $\|(A - z I)^{-1}\|^{-1}$ can be computed in arbitrary norm, when $A$ is close to the diagonal matrix.
\begin{lemma}\label{mA-est}
Let $n\in\N$, $z\in\C$ and $A\in\C^{n\times n}$ be given.
Let $A$ be decomposed into $A=J+E$ where $J$ is a diagonal
matrix and $E$ equals zero on the diagonal. Assume that  $J-z\cdot I_{\C^{n}}$ is invertible and
$\|(J-z\cdot I_{\C^{n}})^{-1}\|^{-1}-\|E\|>0$. Then
\[
\|(A-z\cdot I_{\C^n})^{-1}\|^{-1}\geq \|(J-z\cdot I_{\C^{n}})^{-1}\|^{-1}-\|E\|.
\]
\end{lemma}
\begin{proof}
It is well-known that for an invertible operator $B$ we
have
\[
(B-C)^{-1} = \sum^{\infty}_{n=0}(B^{-1}C)^n B^{-1}\quad\text{for
}\; C\in\C^{n\times n}\;:\;\|C\|<1/\|B^{-1}\|.
\]
Hence, if $\|C\|<1/\|B^{-1}\|$,
then
\begin{eqnarray*}
  \|(B-C)^{-1}\| \leq \frac{\|B^{-1}\|}{1 - \|B^{-1}\| \cdot
  \|C\|},
\end{eqnarray*}
so we obtain
\begin{equation}
  \|(B-C)^{-1}\|^{-1} \geq \frac{1}{\|B^{-1}\|} (1-\|B^{-1}\| \cdot
  \|C\|) = \frac{1}{\|B^{-1}\|} - \|C\|. \label{eq:normA-b-inv}
\end{equation}

From \eqref{eq:normA-b-inv} applied to $B=J-z
I_{\C^{n}}$ and $C=-E$ we get assertion of the lemma.
\end{proof}

Now we present results about the location of the eigenspaces.
\begin{theorem}\label{thm:block-eigenspaces}
Let $k,n\in\N$ such that $0\leq k\leq n$ and
$A\in\LL(\C^k\times\C^{n-k})$ be given in the block from by
\[
A:=
\begin{bmatrix}
 A_{11} & A_{12}\\
 A_{21} & A_{22}
\end{bmatrix},
\]
where $A_{11}\in\LL(\C^k)$, $A_{12}\in\LL(\C^k,\C^{n-k})$, $A_{21}\in\LL(\C^{n-k},\C^{k})$ and $A_{22}\in\LL(\C^{n-k})$. Assume that $A_{22}$ is
invertible, $d:=\|A_{22}^{-1}\|^{-1}-\|A_{11}\|>0$ and
$d^2-4\|A_{12}\|\|A_{21}\|>0$. Then there exists a unique direct
sum decomposition $\C^k\times\C^{n-k}=F_1\oplus F_2$, such that
$F_1$ and $F_2$ are $A$-invariant subspaces $F_1$, $F_2$,  $\dim
F_1=k$, $\dim F_2=n-k$ and
\begin{equation}\label{cor_eq:1}
  \begin{aligned}
F_1 & \subset\left\{(x_1, x_2)\in\C^k\times\C^{n-k} : \|x_2\|\leq\frac{2\|A_{21}\|}{d}\|x_1\| \right\},\\
F_2 & \subset\left\{(x_1, x_2)\in\C^k\times\C^{n-k} :
\frac{2\|A_{21}\|}{d}\|x_1\|\leq\|x_2\| \right\}.
  \end{aligned}
\end{equation}
Moreover, we have
\begin{equation}\label{cor_eq:2}
\sigma(A|_{F_1})\subset\overline{B}\left(0,\|A_{11}\|+\frac{2\|A_{12}\|\cdot\|A_{21}\|}{d}\right),\quad
\sigma(A|_{F_2})\subset\C\setminus
B\left(0,\|A_{22}^{-1}\|^{-1}-\frac{d}{2} \right).
\end{equation}
\end{theorem}
\begin{proof}
Let $c:=\|A_{21}\|$. If $c=0$, the assertion holds. Assume that
$c\neq 0$. By Lemma \ref{lem:2} we get
$A\in\D_{d/(2c)}(\C^k\times\C^{n-k})$, and from Theorem
\ref{thm:cone-main-eigen-location} we know that  exists a direct sum decomposition
$\C^k\times\C^{n-k}=F_1\oplus F_2$ such that $\dim F_1=k$, $\dim
F_2=n-k$ and $F_1$, $F_2$ are invariant. The properties
\eqref{cor_eq:1} and \eqref{cor_eq:2} are  consequences of Theorem
\ref{thm:cone-main-eigen-location} and  Theorem \ref{thm:formulas-contr-exp} and Definition \ref{def:1},
respectively.
\end{proof}

\begin{corollary}\label{cor:block-eigenspaces}
We use the same notation and decomposition of the matrix A as in Theorem~\ref{thm:block-eigenspaces}.
Assume that for some $z\in \C$ matrices $A_{11} - zI_{\C^k}$, $A_{22}-
zI_{\C^{n-k}}$ are invertible and $d:=\|(A_{22}-
zI_{\C^{n-k}})^{-1}\|^{-1}-\|A_{11}-zI_{\C^{k}}\|>0$,
$d^2-4\|A_{12}\|\|A_{21}\|>0$. Then there exists a unique direct
sum decomposition $\C^k\times\C^{n-k}=F_1\oplus F_2$ into
$A$-invariant subspaces $F_1$, $F_2$ such that $\dim F_1=k$, $\dim
F_2=n-k$ and
\begin{align*}
F_1 & \subset\left\{(x_1, x_2)\in\C^k\times\C^{n-k} : \|x_2\|\leq\frac{2\|A_{21}\|}{d}\|x_1\| \right\},\\
F_2 & \subset\left\{(x_1, x_2)\in\C^k\times\C^{n-k} :
\frac{2\|A_{21}\|}{d}\|x_1\|\leq\|x_2\| \right\}.
\end{align*}
Moreover, we have
\begin{align*}
\sigma(A|_{F_1}) & \subset\overline{B}\left(z,\|A_{11}-zI_{\C^{k}}\|+\frac{2\|A_{12}\|\cdot\|A_{21}\|}{d}\right), \\
\sigma(A|_{F_2}) & \subset\C\setminus B\left(z,\|(A_{22}-
zI_{\C^{n-k}})^{-1}\|^{-1}-\frac{d}{2} \right).
\end{align*}
\end{corollary}

\subsection{Gerschgorin theorems}
\label{subsec:Ger-th}

For to the convenience of the reader, in this section we recall the Gerschgorin theorem and its modifications.

We have a matrix $A$ which has a block structure
\[
  A = \left[
        \begin{array}{cccc}
          A_{11} & A_{12} & \cdots & A_{1n} \\
          A_{21} & A_{22} & \cdots & A_{2n} \\
          \cdots & \dots & \dots & \dots \\
          A_{n1} & A_{n2} & \dots & A_{nn} \\
        \end{array}
      \right],
\]
where $A_{ij}$ are  matrices and $A_{ii}$ are square matrices.

Let $V=\bigoplus_{i=1}^n V_i $,
where $V_i$ are finite dimensional vector spaces over
$\mathbb{C}$, and $A:V \to V$ be decomposed into blocks
$A_{ij}:V_j \to V_i$ $i,j=1,2,\dots,n$, so that for $v=v_1 + \dots
+ v_n$, where $v_i \in V_i$ holds
\begin{equation}
  A(v_1 + \dots + v_n)= \sum_i \sum_j A_{ij}v_j.
\end{equation}

We define Gerschgorin disks $G_i(A)$ for the block matrix $A$ by
\begin{eqnarray*}
  R_i(A)&=&  \sum_{j,j\neq i} \|A_{ij}\| , \\
  G_i(A)&=&\{\lambda \in \mathbb{C} \ : \ \mbox{$A_{ii}-\lambda I_{i}$ exists and } \|(A_{ii}-\lambda I_{i})^{-1}\|^{-1} \leq R_i(A)\}, \quad i=1,\dots,n,
\end{eqnarray*}
where $I_{V_i}$ is an identity map on $V_i$.  If $A$ is known from the context, then we will usually drop $A$ and write just $R_i$ and $G_i$.
Similarly, we write $I$ instead of $I_{V_i}$.

Theorem below we present the generalizations of Gershgorin Theorems due to Feingold and Varga \cite{Gen}.
\begin{theorem} \cite[Theorem 2]{Gen}
\label{thm:gen-gersch-thm1}
  \[
     \sigma(A) \subset \bigcup_{i=1}^n G_i.
  \]
\end{theorem}

\begin{theorem}  \cite[Theorem 4]{Gen}
\label{thm:gersz-eigenval}
Assume that $J \subset \{1,\dots,n\}$ is such that
\[
  \left(\bigcup_{j \in J} G_j \right) \cap  \left(\bigcup_{j \notin J} G_j
  \right)= \emptyset.
\]
Then the number of eigenvalues of $A$ (counting with
multiplicities) contained in $\left(\bigcup\limits_{j \in J} G_j \right)$ is
equal to $\sum\limits_{j \in J} \mbox{dim}\, V_j$.
\end{theorem}

Now we give a theorem about the location of the eigenvectors based on the   Wilkinson  argument~\cite{Wilk}.
\begin{theorem}\label{thm:gen-gersz-eigenvect-estm}
 Assume that for some $j\in\{1,\ldots,n\}$
\[
 G_j \cap G_k =\emptyset, \quad \mbox{for $k=1,2,\dots,n$, $k \neq
 j$}.
\]

Then if $v=(v_1+\dots+v_n)$ is an eigenvector corresponding to
$\lambda \in G_j$, then $\|v_k\| \leq \|v_j\|$ for $k=1,\dots,n$.
\end{theorem}
\begin{proof}
To show that $\|v_k\| \leq \|v_j\|$ we will reason by the
contradiction. Assume that for some $i \neq 0$ holds $\|v_i\|\geq
\|v_k\|$ for $k=1,\dots,n$ and $\|v_i\| > \|v_j\|$. We will apply
the basic argument from the generalized Gerschgorin theorem
(Theorem~\ref{thm:gen-gersch-thm1}) to prove that $\lambda \in
G_i$. This will lead to a contradiction, because $\lambda \in
G_j$, hence $\lambda \in G_j \cap G_i \neq \emptyset$.

We have
\begin{eqnarray*}
  \lambda v_i&=&A_{ii} v_i + \sum_{k \neq i} A_{ik} v_k \\
  (\lambda I - A_{ii})v_i &=& \sum_{k \neq i} A_{ik} v_k  \\
  \|(\lambda I - A_{ii})^{-1}\|^{-1} \|v_i\| &\leq&  \sum_{k \neq i} \|A_{ik}\| \|v_k\| \\
   \|(\lambda I - A_{ii})^{-1}\|^{-1}   &\leq&  \sum_{k \neq i} \|A_{ik}\|
   \frac{\|v_k\|}{\|v_i\|} \leq \sum_{k \neq i} \|A_{ik}\|
\end{eqnarray*}
hence $\lambda \in G_i$. We obtained the contradiction.  This
finishes the proof.
\end{proof}

One of the easiest ways to improve the estimation of the eigenvalues from the Gerschgorin theorem is through the scaling the basis of our domain.
This approach is well known and can be found
in the original article of Gerschgorin \cite{G}.

Assume, that we have matrix $A\in\C^{n\times n}$ and let
$\x{x}=(x_1,\ldots,x_n)^T\in\R^n$ such that $x_i>0$ for all
$i\in\{1,\ldots, n\}$. With this vector $\x{x}$ we define the
matrix $X\in\R^{n\times n}$ with the elements of $\x{x}$ on the
leading diagonal, and $0$ elsewhere. Note, that the matrix $X$ is
nonsingular and matrix $X^{-1}AX$ is similar to $A$ therefore
$\sigma(X^{-1}AX)=\sigma(A)$. If $A=[a_{ij}]_{1\leq i,j\leq n}$,
then
\[
X^{-1}AX=\left[\frac{a_{ij}x_j}{x_i}\right]_{1\leq i,j\leq
n}\]
and
\[
G_i=\overline{B}\Big(a_{ii},\sum\limits_{j \neq i}
\frac{|a_{ij}|x_j}{x_i}\Big) \quad\mbox{ for $i=1,\dots,n$}.
\]

\subsection{Example}
In the following example we consider a matrix with  multi-dimensional block for which we estimate eigenspaces.

\begin{example}\label{ex:subspace}
Consider the matrix $A\in\LL(\C^2\times\C^2)$ be given by
\[
  A=
\begin{bmatrix}
A_{11} & A_{12} \\
A_{21} & A_{22}
\end{bmatrix}
= \left[
\begin{array}{cc|cc}
 0. & 0.15 & 0.11 & 0.02 \\
 0.2 & 0. & 0.1 & 0.05 \\ \hline
 0.01 & 0.025 & 0. & 1.5 \\
 0.15 & 0.05 & 1. & 0. \\
\end{array}
\right].
\]
We have $\|A_{11}\|_{\infty}=0.2$, $\|A_{12}\|_{\infty}=0.15$, $\|A_{21}\|_{\infty}=0.2$.
From Theorem~\ref{thm:block-eigenspaces}
($d=\|(A_{22}^{-1}\|_{\infty}^{-1}-\|A_{11}\|_{\infty}=1-0.2=0.8>0$ and
$d^2-4\|A_{12}\|_{\infty}\|A_{21}\|_{\infty}=0.52>0$) we know that there exist
eigenspaces $F_1$ and $F_2$, which satisfy
\begin{align*}
F_1 & \subset\left\{(x_1, x_2)\in\C^2\times\C^{2} : \|x_2\|\leq 0.5 \|x_1\| \right\},\\
F_2 & \subset\left\{(x_1, x_2)\in\C^2\times\C^{2} : \|x_1\|\leq
2\|x_2\|\right\}.
\end{align*}
 and
$\sigma(A_{F_1})\subset\overline{B}(0,0.275)$,
$\sigma(A_{F_2})\subset\C\setminus B(0,0.6)$ (see Figure
\ref{fig2:b}).

 \begin{figure}[H]
\setcounter{subfigure}{0}
  \centering
 \subfloat[Gerschgorin circles.]{\includegraphics[width=0.48\textwidth]{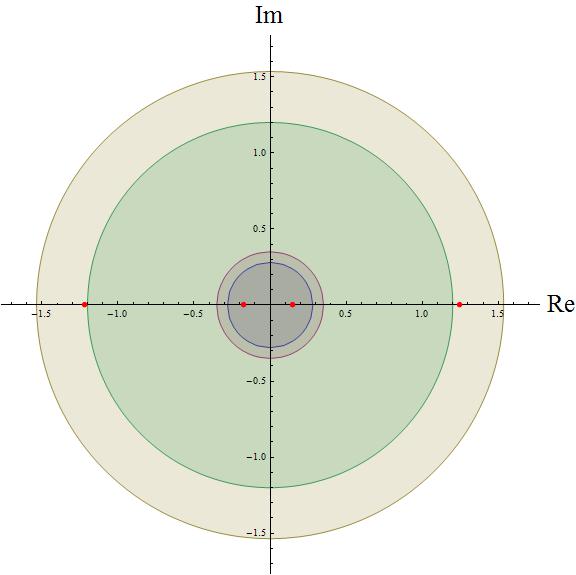}}
 \quad
 \subfloat[Estimates based on Theorem~\ref{thm:block-eigenspaces}. The white annulus in does not contain any eigenvalue.\label{fig2:b}]{\includegraphics[width=0.48\textwidth]{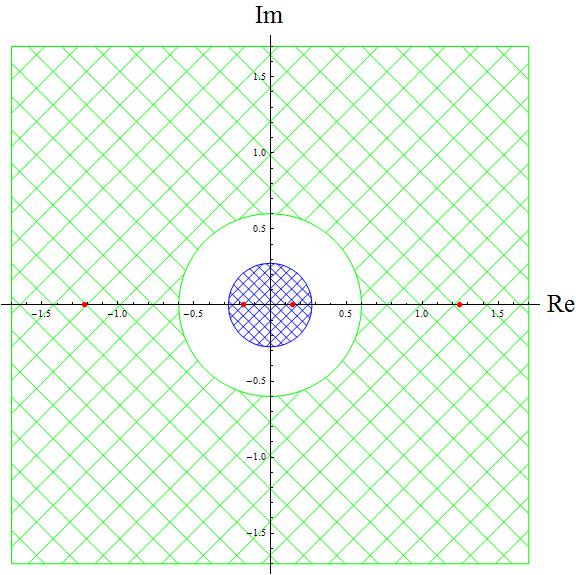}}
\caption{Gerschgorin and our circles with approximate eigenvalues in Example \ref{ex:subspace}.}
 \end{figure}

Observe that when using the Gerschgorin theorem with one-dimensional blocks with scalings, as described at the end of Section~\ref{subsec:Ger-th}, we will not be able to separate the spectrum of $A$, because the centers of  Gerschgorin
circles are located at zero.

\medskip

Now we discuss what happens when we use the generalized Gerschgorin theorems from \cite{Gen}. First rescale the matrix $A$ by
\(
X=
\begin{bmatrix}
2 & 0 \\ 0 & I
\end{bmatrix}
\)
(we take the same rescaling as in our method, see Remark \ref{rem:good-eps}) to get
\[
\tilde{A}=X^{-1}AX=
\begin{bmatrix}
A_{11} & \frac{1}{2}A_{12} \\ 2A_{21} & A_{22}
\end{bmatrix}.
\]
We use the Theorems~\ref{thm:gen-gersch-thm1} and~\ref{thm:gersz-eigenval} applied to the above block decomposition, and
obtain the generalized Gerschgorin disks:
\begin{eqnarray*}
  G_1(\tilde{A})&=& \left\{\lambda\in\C : \|(A_{11} - \lambda I)^{-1}\|_{\infty}^{-1} \leq \frac{1}{2}\|A_{12}\|_\infty\right\}, \\
  G_2(\tilde{A})&=&\left\{\lambda\in \C : \|(A_{22} - \lambda I)^{-1}\|_{\infty}^{-1} \leq 2\|A_{21}\|_\infty\right\}.
\end{eqnarray*}
We want to show that $G_1(\tilde{A}) \cap G_2(\tilde{A}) = \emptyset$. Let us we check that $G_1(\tilde{A})\subset\overline{B}(0,0.25)$.
We have
\[
(A_{11}-\lambda I)^{-1}=\frac{1}{\lambda^2-0.03}
\begin{bmatrix}
-\lambda & -0.15 \\
-0.2 & -\lambda
\end{bmatrix},
\]
so we get
\[
\|(A_{11}-\lambda I)^{-1}\|_{\infty}^{-1}=\frac{|\lambda^2-0.03|}{0.2+|\lambda|}.
\]
For $\lambda\in G_1(\tilde{A})\subset\C$ we have
\[
\frac{|\lambda^2-0.03|}{0.2+|\lambda|}\leq 0.075.
\]
Performing simple mathematical operations and changing the coordinate system to the polar one we obtain
\[
40000 r^4 - 15 r (160 r \cos (2 q)+15 r+6) +27\leq 0, \quad r=|\lambda|\in[0,\infty),\ \varphi\in[0,2\pi).
\]
Solving the above inequality we get
\[
\sup r = \frac{3}{80} \left(1+\sqrt{33}\right)< \frac{21}{80}.
\]
This means that $G_1(\tilde{A})\subset\overline{B}(0,21/80)$. Now we show that $\lambda\notin G_2(\tilde{A})$ for an arbitrary $\lambda\in\overline{B}(0,21/80)$.

Indeed we have
\[
(A_{22}-\lambda I)^{-1}=\frac{1}{\lambda^2-1.5}
\begin{bmatrix}
-\lambda & -1.5 \\
-1 & -\lambda
\end{bmatrix}.
\]
It is easy to see that for $\lambda\in\overline{B}(0,21/80)$ we have
\[
\|(A_{22}-\lambda I)^{-1}\|_{\infty}< \frac{\frac{3}{2}+\frac{21}{80}}{\frac{3}{2}-\left(\frac{21}{80}\right)^2}=\frac{3760}{3053}.
\]
Hence
\[
\|(A_{22}-\lambda I)^{-1}\|_{\infty}^{-1}>\frac{3053}{3760}>\frac{8}{10}, \qquad \lambda\in G_1(\tilde{A})\subset \overline{B}(0,21/80).
\]
Finally, we get $G_1(\tilde{A}) \cap G_2(\tilde{A}) = \emptyset$ (see Figure \ref{fig2:c}) and therefore we obtain from Theorem~\ref{thm:gen-gersch-thm1} and~\ref{thm:gersz-eigenval} that two eigenvalues belong to $G_1(\tilde{A})$ while the remaining two eigenvalues are inside $G_2(\tilde{A})$.
As one can see, we get better estimation for eigenvalues close to $0$ from generalized Gerschgorin theorem with scaling $r=2$, than from Theorem~\ref{thm:block-eigenspaces} but by generalized Gerschgorin theorem we can not get eigenspaces.

\begin{figure}[H]
  \centering
      \includegraphics[width=0.55\textwidth]{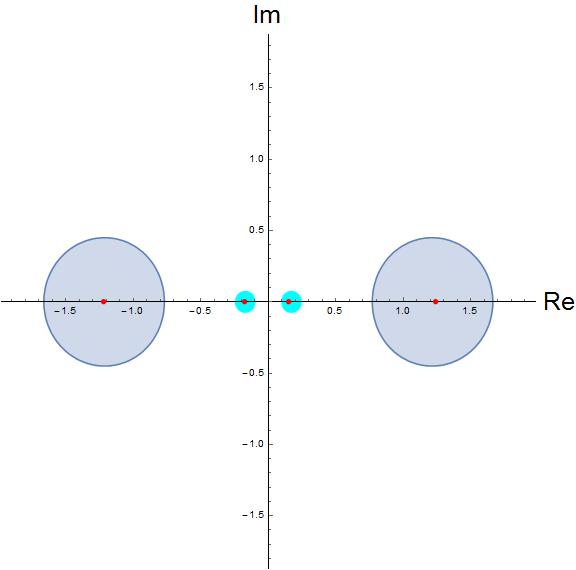}
  \caption{Generalized Gerschgorin circles: $G_1(\tilde{A})$ -- greater circles and $G_2(\tilde{A})$ -- smaller ones in Example \ref{ex:subspace} (compare Fig.~\ref{fig2:b}).}
  \label{fig2:c}
\end{figure}
\end{example}

\section{Comparisons in the case of the isolated Gerschgorin disk}
\label{sec:comparison}

In this section we compare our method of cones with the Gerschgorin theorem with rescaling of the basis, when
trying to estimate an eigenvalue in an isolated Gerschgorin disk  and corresponding eigenvector. Throughout this section we will use the $\|\cdot\|_\infty$ norm.

\subsection{The isolation of first Gerschgorin disk implies that the matrix $A-a_{11} I$ is dominating}
When applying Theorem~\ref{thm:gersz-eigenval} with  the splitting $\C \oplus \C^{n-1}$  we will have two generalized Gerschgorin disks
\begin{eqnarray*}
  G_1(A)&=&\overline{B}(a_{11},\|A_{12}\|_\infty)=\overline{B}(a_{11},\sum_{j\neq 1}|a_{1j}|), \\
  G_2(A)&=&\{ \lambda \in \mathbb{C} \ : \   \|(A_{22} - \lambda I)^{-1}\|_{\infty}^{-1} \leq \max_{j=2,\dots,n} |a_{j1}| \}.
\end{eqnarray*}

Now we develop computable bounds for $G_2(A)$.
\begin{lemma}
\label{lem:mA-Gersz}
Let $A=[a_{ij}]\in\C^{n\times n}$.
Then
\begin{equation}
  \min_i (|a_{ii}| - \sum_{j \neq i} |a_{ij}|) \leq \sup\{ \lambda \in \mathbb{R}\ | \ \forall x\in \C^n \quad \|Ax\|_\infty  \geq \lambda \|x\|_\infty \}.
    \label{eq:m(A)}
\end{equation}
If $A$ is invertible, then
\[
 \min_i (|a_{ii}| - \sum_{j \neq i} |a_{ij}|) \leq \frac{1}{\|A^{-1}\|_\infty}.
\]
\end{lemma}
\begin{proof}
Let
\[
S:= \min_i (|a_{ii}| - \sum_{j \neq i} |a_{ij}|).
\]
 Let us take any $x \in \mathbb{C}^n$, such that $\|x\|=1$. Let $i$ be such that $|x_i|=1$.   We have
\[
   |(Ax)_i| \geq (|a_{ii}| |x_i| - \sum_{j \neq i} |a_{ij}| \cdot |x_j|) \geq (|a_{ii}| - \sum_{j \neq i} |a_{ij}|) \geq  S >0.
\]
 Hence
\[
   \|Ax\|_\infty \geq S.
\]
This establishes \eqref{eq:m(A)}.

For the second part observe that, if $A$ is invertible, then
\[
\sup\{ \lambda \in \mathbb{R}\ | \ \forall x\in \C^n \quad \|Ax\|_\infty  \geq \lambda \|x\|_\infty \} = \frac{1}{\|A^{-1}\|_\infty}.
\]
\end{proof}

From Lemma~\ref{lem:mA-Gersz} it follows that
 \begin{eqnarray*}
  G_2(A) &\subset&  \{ \lambda \in \mathbb{C} \ | \   \min_{i=2,\dots,n}( |a_{ii} - \lambda| - \sum_{j \neq 1,i} |a_{ij}|) \leq \max_{j=2,\dots,n} |a_{j1}| \} \\
   &=& \{ \lambda \in \mathbb{C} \ | \   \exists i=2,\dots,n \   |a_{ii} - \lambda| \leq \sum_{j \neq 1,i} |a_{ij}| + \max_{j=2,\dots,n} |a_{j1}| \}.
 \end{eqnarray*}
So we see that $G_1(A) \cap G_2(A) =\emptyset $ if the following condition holds
\begin{equation}
  |a_{11}-a_{ii}| > \sum_{j\neq 1}|a_{1j}| + \sum_{j \neq 1,i} |a_{ij}| + \max_{j=2,\dots,n} |a_{j1}| \quad\mbox{ for all $i=2,\dots,n$.} \label{eq:GenGerszEstm}
\end{equation}

If we will use the classical Gerschgorin theorem, i.e.  blocks are one-dimensional, then to have $G_1 \cap G_i=\emptyset$ for
\begin{equation}
  |a_{11}-a_{ii}| > R_1 + R_i = \sum_{j \neq 1} |a_{1j}|  + \sum_{j \neq i} |a_{ij}| \quad\mbox{ for all $i=2,\dots,n$.} \label{eq:StandardGerszEstm}
\end{equation}
Observe that in both cases we have the same Gerschgorin disk $G_1$, so the bound for the first eigenvalue will be the same, provided we have
empty intersections with other disks.
Observe that (\ref{eq:GenGerszEstm}) implies (\ref{eq:StandardGerszEstm}).

Now we show one of the main results of this paper, which states that if a matrix $A=[a_{ij}]$ has an isolated Gerschgorin disk $G_1$, then $A-a_{11}I$ is dominating (relative to the splitting $\C\times\C^{n-1}$) and under very mild assumptions the bound obtained from the method of cones is better that the
one from the Gershgorin theorem.

\begin{theorem}
\label{thm:gerssh-dominating}
Let $A\in\C^{n\times n}$ be given by the formula
\begin{equation*}
A=
\left[
\begin{array}{c|ccc}
a_{11} & a_{12} & \ldots & a_{1n} \\ \hline
a_{21} & a_{22} & \ldots & a_{2n} \\
\vdots & \vdots & \ddots & \vdots \\
a_{n1} & a_{n2} & \ldots & a_{nn}
\end{array}
\right]
=
\begin{bmatrix}
A_{11} & A_{12}  \\
A_{21} & A_{22}
\end{bmatrix}.
\end{equation*}
Assume that the matrix $A$ satisfies inequality \eqref{eq:StandardGerszEstm}.
Then the matrix $A-a_{11}I$ is dominating. Moreover, we have
\[
\co{A-a_{11}I}\ \leq \sum\limits_{j\neq 1}|a_{1j}| < \ex{A-a_{11}I}.
\]
and if $\sum\limits_{j\neq 1}|a_{1j}| >0$, then
\[
\co{A-a_{11}I}\ < \sum\limits_{j\neq 1}|a_{1j}|.
\]
\end{theorem}

\begin{proof}
Without any loss of the generality we can assume that $a_{11}=0$. Let us denote $V=\mathbb{C} \oplus \mathbb{C}^{n-1}$

In order to estimate $\co{A}$ we will first find a bound for the set $Z$ of $\x{x}$, such that $A\x{x} \in \co{V}$. Then
we will compute $\co{A}$ on $Z$.

Let
\begin{equation}
\label{eq:thm-better-1}
\delta_k=|a_{kk}| - \sum_{j \neq k}|a_{kj}| - \sum_{j \neq 1}|a_{1j}|.
\end{equation}
From our assumptions it follows that
\[
 \delta=\min_{k=2,\dots,n} \delta_k >0.
\]
Let $\epsilon>0$ be such that
\begin{equation}
\label{eq:thm-better-2}
\epsilon |a_{k1}| < \delta_k, \quad k=2,\dots,n.
\end{equation}
Assume now that $\x{x}=(x_1,\x{x}_2)$ such that $|x_1| \leq (1+\epsilon) \|\x{x}_2\|_\infty$. We will show that $Ax \notin \co{V}$.

We can assume that $\|x_2\|_\infty=1$ and $|x_k|=1$. Then we have
\begin{align*}
  |(Ax)_k| & \geq |a_{kk}| - \sum_{j \notin \{1,k\}} |a_{kj}| - (1+\epsilon) |a_{k1}|= |a_{kk}| - \sum_{j \neq k} |a_{kj}| - \epsilon |a_{k1}| \\[0.3em]
& \overset{by\,\eqref{eq:thm-better-1}}{=} \sum_{j \neq 1} |a_{1j}| + \delta_k - \epsilon |a_{k1}|
\overset{by\,\eqref{eq:thm-better-2}}{>}  \sum_{j \neq 1} |a_{1j}| \geq \|A_{12} x_2\|=\|(A\x{x})_1\|_\infty.
\end{align*}
Hence $A\x{x} \notin \co{V}$.

Therefore, if $A\x{x} \in \co{V}$, then $|x_1| > (1+\epsilon)\|\x{x}_2\|_\infty$. In particular, we obtain
\begin{equation}
  \mbox{if $A\x{x} \in \co{V}$, then} \quad \norm{\x{x}}=|x_1| > (1+\epsilon)  \|\x{x}_2\|_\infty. \label{eq:coA-set-estm}
\end{equation}

Now we are ready to estimate $\co{A}$. Let  $\x{x}=(x_1, \x{x}_2)$ is such that $A\x{x} \in \co{V}$, then
\begin{align*}
\norm{A\x{x}} & = \|A_{12}\x{x}_2\|_{\infty}\leq \|A_{12}\|_{\infty}\cdot \|\x{x}_2\|_{\infty}\\[0.3em]
& \overset{by\,\eqref{eq:coA-set-estm}}{\leq}
 \|A_{12}\|_{\infty} \frac{\norm{\x{x}}}{1+\epsilon} =  \frac{1}{{1+\epsilon}}\sum\limits_{j\neq 1}|a_{1j}| \norm{\x{x}}.
\end{align*}
Hence $\co{A}\ \leq \sum\limits_{j\neq 1}|a_{1j}|$, but if $\sum\limits_{j\neq 1}|a_{1j}| >0$, then $\co{A}\ < \sum\limits_{j\neq 1}|a_{1j}|$.

Now we estimate $\ex{A}$. We will use Lemma~\ref{lem:coA-exA-alter}.
Let's take arbitrary $\x{x}=(x_1,\x{x}_2)$ such that $\|\x{x}_2\|_{\infty}=1$ and $|x_1|\leq 1$. Let $k=2,\ldots, n$ be such that $|x_k|=1$. From \eqref{eq:StandardGerszEstm} we obtain
\[
|(Ax)_k|\geq |a_{kk}|-\sum\limits_{j\neq k}|a_{kj}|  \overset{by\,\eqref{eq:StandardGerszEstm}}{>} \sum\limits_{j\neq 1}|a_{1j}|.
\]
Hence $\|Ax\|_{\infty}>\sum\limits_{j\neq 1}|a_{1j}|$.
Therefore we have shown that
\[
\ex{A}>\sum\limits_{j\neq 1}|a_{1j}|.
\]
\end{proof}

\begin{remark}
Observe that from Theorem~\ref{thm:gerssh-dominating} we know that our method of cones (i.e. Theorem~\ref{thm:cone-main-eigen-location}) can be used for all matrices which have an isolated Gerschgorin disk. Moreover, we obtain
\[
|\lambda-a_{11}|\leq\ \co{A-a_{11}I}\ \leq  \frac{1}{1+\epsilon} R_1 = \frac{1}{1+\epsilon} \sum\limits_{j\neq 1}|a_{1j}|.
\]
This means that, if $R_1$ (the radius of the first Gerschogorin disk) is nonzero, then the estimate of the first eigenvalue from our method based on cones is better than the one obtained from the Gerschgorin theorem. This is  also
valid for all possible rescalings in the application of the Gerschogorin theorem, we should apply the same scaling in the method of cones.
\end{remark}

\subsection{Comparison of Theorem~\ref{thm:cone-single-eigenval} with the Gerschgorin theorems}

In the proof of Theorem~\ref{thm:cone-single-eigenval} we  applied Theorem~\ref{thm:cone-main-eigen-location} to the matrix $A - a_{11}I$  to estimate the size of the eigenvalue, $\lambda_1$, close to $0$. We looked for possibly large parameter $r$, such that $A-a_{11}I$ is $r$-dominating and  then we obtain
\[
   |a_{11} - \lambda_1| \leq \co[r]{A} \leq \frac{\|A_{12}\|}{r}.
\]
 This is exactly $G_1$ obtained from the Gerschgorin theorem for $\tilde{A}_r$.

  The  optimization with respect of $r$ performed in the proof  of the Theorem~\ref{thm:cone-single-eigenval}
 to obtain the formula  can be also repeated by suitable rescaling using the original Gerschgorin theorem as long $G_1(\tilde{A}_r)$ is disjoint from
  other Gerschgorin disks for $\tilde{A}_r$.  Therefore both approaches differ only with the range of $r$'s over which the optimization can be performed.
  In fact we are only interested in the upper bound for $r$ in both methods.

Let $(1,\delta_2,\dots,\delta_n)$ be the eigenvector corresponding to $\lambda_1$.  We obtain from  Theorem~\ref{thm:r-dom-spectrum-gap}
 the bound $1 \geq r \|(\delta_2,\dots,\delta_n\|$, while from Theorem~\ref{thm:gen-gersz-eigenvect-estm} applied to $\tilde{A}$ after returning to the original base we have $|\delta_i| \leq 1/r$. Hence the result is the same for the method based on cones and the Gerschgorin Theorem.

\medskip

The example below demonstrates that it is possible to use the Gerschgorin theorem to isolate and estimate the eigenvector and eigenvalue, while
assumptions of Theorem~\ref{thm:cone-single-eigenval} and also assumptions of the generalized Gerschgorin Theorems~\ref{thm:gen-gersch-thm1} and~\ref{thm:gersz-eigenval} are not satisfied. This appears to contradict Theorem~\ref{thm:gerssh-dominating}, but it does not, because
in the proof Theorem~\ref{thm:cone-single-eigenval} we used an expression for $\co{A}$ from Lemma~\ref{thm:formulas-contr-exp}, which turns out to be
an overestimation, see also Example~\ref{ex:dominat-exactly}.
By Theorem \ref{thm:gerssh-dominating} we know that the considered matrix  is dominating, hence we can estimate the eigenpair using Theorem~\ref{thm:cone-main-eigen-location}, see Example \ref{ex:not-G-better}.

\begin{example}
\label{ex:G-better}
Let $A\in\LL(\C\times\C^2)$ be given by the formula
\[
A=
\begin{bmatrix}
A_{11} & A_{12} \\
A_{21} & A_{22}
\end{bmatrix}
= \left[
\begin{array}{c|cc}
0 & 1 & 0 \\ \hline
0.5 & 2 & 0 \\
50 & 0 & 100
\end{array}
\right].
\]
The classical Gerschgorin disks are
\[
  G_1=\overline{B}(0,1), \quad G_2=\overline{B}(2,0.5), \quad G_3=\overline{B}(100,50).
\]
It is clear that  they are mutually disjoint, hence from the Gerschgorin theorem there exists an eigenvalue $\lambda$, $|\lambda|\leq 1$.

Now, we look at our Theorem~\ref{thm:cone-single-eigenval} to estimate the eigenvalue close to $0$. We have $A_{11}=0$ and
\begin{eqnarray*}
  \|A_{12}\|_\infty=1, \quad \|A_{21}\|_{\infty}=50, \quad \|(A_{22} - A_{11}\cdot I)^{-1}\|_\infty=0.5, \\
   \|(A_{22} - A_{11}\cdot I)^{-1}\|_\infty^{-2} - 4 \|A_{12}\| \cdot \|A_{21}\|_\infty= 4 - 200 < 0.
\end{eqnarray*}
Therefore assumptions of  Theorem~\ref{thm:cone-single-eigenval} are not satisfied.

Observe that also assumptions of the generalized Gerschgorin Theorems~\ref{thm:gen-gersch-thm1} and~\ref{thm:gersz-eigenval}  for the decomposition given above are not satisfied. Our generalized Gerschgorin disks are
\begin{eqnarray*}
  G_1(A)&=&\overline{B}(0,1), \\
  G_2(A)&=&\{\lambda \ :\ \|(A_{22} - \lambda I)^{-1}\|_{\infty}^{-1} \leq 50\}.
\end{eqnarray*}
We have
\begin{equation*}
  (A_{22} - \lambda I)^{-1}= \frac{1}{(100-\lambda)(2-\lambda)}
\begin{bmatrix}
100-\lambda & 0 \\
0 & 2-\lambda\\
\end{bmatrix},
\end{equation*}
hence
\begin{equation*}
\|(A_{22} - \lambda I)^{-1}\|_{\infty}^{-1}=\min(|2-\lambda|,|100-\lambda|).
\end{equation*}
It is easy to see that $G_1(A)\cap G_2(A)\neq\emptyset$, therefore we cannot use these theorems.

\medskip

\textbf{Rescaling:} When applying our method based on cones we should look for the largest $r$ such that $\tilde{A}_r$ is $1$-dominating and when using the Gerschgorin theorem we look for $r$
such that $G_1(\tilde{A}_r)$ have empty intersection with others Gerschgorin circles for $\tilde{A}_r$.

For the Gerschorin disks we need to have the following inequalities
\[
  \frac{1}{r} < 2 - r/2, \quad \frac{1}{r} < 100 - 50r.
\]
We obtain $\sup r =1 + \sqrt{\frac{49}{50}} \approx 2$. Hence we obtain bound $|\lambda| \leq \approx 1/2$.

For the approach based on cones we need to find largest $r$, such that $\tilde{A}_r$ is $1$-dominating. Using Theorem~\ref{thm:formulas-contr-exp} we
obtain the following condition
\[
\frac{1}{r}  \|A_{12}\|= \frac{1}{r} < \|A_{22}^{-1}\|^{-1} - r \|A_{21}\| = 2 - 50r.
\]
Easy computations show that no such $r$ exists in this case.
Similar effect we get if we use the generalized Gerschgorin theorem.
\end{example}

In the following example we show that despite the fact that the matrix $A$ from Example~\ref{ex:G-better} does not satisfy the assumptions of Theorem~\ref{thm:cone-single-eigenval}, we can use our method of cones (we apply Theorem~\ref{thm:cone-main-eigen-location}) to estimate the eigenvalue close to zero.

\begin{example}
\label{ex:not-G-better}
Recall that $A\in\LL(\C\times\C^2)$ of Example~\ref{ex:G-better} was given by the formula
\[
A=
\left[
\begin{array}{c|cc}
0 & 1 & 0 \\ \hline
0.5 & 2 & 0 \\
50 & 0 & 100
\end{array}
\right].
\]
From Theorem \ref{thm:gerssh-dominating} we know that the matrix $A$ is dominating, so we can estimate the eigenvalue $\lambda$ close to zero by $|\lambda|\leq\co{A}$ (see Theorem~\ref{thm:cone-main-eigen-location}).
From Lemma~\ref{lem:coA-exA-alter} we have
\begin{equation*}
\co[]{A} = \dfrac{1}{\min\left(\|x\|_{\infty} \ \mbox{ for } \x{x}\in\R^3\ \mbox{ such that }\  \|A\x{x}\|_{>1}\leq  \|A\x{x}\|_{\leq 1}=1 \right)},
\end{equation*}
where $\|\x{x}\|_{\leq k}:=\max\limits_{i\leq
k}|x_i|$ and $\|\x{x}\|_{>k}:=\max\limits_{i> k}|x_i|$ for $\x{x}=(x_1,\ldots, x_k,\ldots, x_n)\in\R^n$, see \eqref{eq:coA-alter}.

The problem to calculate this constant comes down to solve simple optimization problem. We obtain \[
\min\left(\|x\|_{\infty} \ \mbox{ for } \x{x}\in\R^3\ \mbox{ such that }\  \|A\x{x}\|_{>1}\leq  \|A\x{x}\|_{\leq 1}=1 \right)=2.
\]
This minimum is realized in the points $\left(-2,1,\frac{99}{100}\right)^T$, $\left(-2,1,\frac{101}{100}\right)^T$, $\left(2,-1,-\frac{101}{100}\right)^T$ and $\left(2,-1,-\frac{99}{100}\right)^T$. Hence $\co[]{A}=\frac{1}{2}$. By  Theorem \ref{thm:cone-main-eigen-location} we get that eigenvalue close to zero satisfies
\[
|\lambda|\leq\co{A}=\frac{1}{2}.
\]

The bound $|\lambda|\leq\frac{1}{2}$ can be obtained also from  the Gerschgorin theorem, see Example \ref{ex:G-better} ({\em 'Rescaling'}). Note that so far we did not improve the matrix $A$ through the scaling
\(
X=
\begin{bmatrix}
r & 0 \\
0 & I
\end{bmatrix}
\)
for $r\in(0,\infty)$.
From Theorem \ref{thm:gerssh-dominating} and again  calculations from Example \ref{ex:G-better} ({\em 'Rescaling'}) we know that our method work even if we rescale our matrix $A$ by the matrix
$X$ for $r<1+\sqrt{\frac{49}{50}}$. For $r=\frac{9}{5}$ we obtain $|\lambda|\leq\co{A}=\frac{9}{26}<\frac{1}{2}$.
\end{example}

In the following two examples in view of the complicated mathematical calculations we will not try to apply the generalized Gerschgorin theorem (in both examples assumptions of Theorems~\ref{thm:gen-gersch-thm1} and~\ref{thm:gen-gersz-eigenvect-estm} are satisfied).
In the first example we construct a matrix such that the matrix $A - a_{11} I$ will be $1$-dominating, while there will be no isolation of the
first Gerschgorin disk.

\begin{example}
\label{ex:mA-better-G}
Let $A\in\LL(\C\times\C^2)$ be given by the formula
\[
A=
\begin{bmatrix}
A_{11} & A_{12} \\
A_{21} & A_{22}
\end{bmatrix}
= \left[
\begin{array}{c|cc}
0 & 0.75 & 0 \\ \hline
\epsilon_1 & 1& 0.5 \\
\epsilon_2 & 0.5 & 100
\end{array}
\right],
\]
where $\epsilon_1$, $\epsilon_2$ are sufficiently small.
Observe that $G_1(A) \cap G_2(A) = \overline{B}(0,0.75)\cap\overline{B}(1,0.5+\epsilon_1) \neq \emptyset$, hence the Gerschgorin theorem does not give us that $\lambda_1 \in G_1(A)$.

It is easy to see that  $A - a_{11} I$ will be $1$-dominating. Indeed  from Theorem~\ref{thm:formulas-contr-exp} we have
\begin{eqnarray*}
  \co[1]{A} \leq \|A_{12}\|= 3/4, \quad \ex[1]{A}  \geq \|A_{22}^{-1}\|^{-1} - \|A_{21}\| \approx 1.
\end{eqnarray*}
Hence $A$ is $1$-dominating and Theorem~\ref{thm:cone-main-eigen-location} implies that $\lambda_1 \in G_1(A)$.

\medskip

\textbf{Rescaling:} We set $\epsilon_1=\epsilon_2=0.1$.
 We optimize by rescaling by $r$. The Gerschgorin disks approach leads to the following inequalities
\[
 \frac{3}{4 r} < 0.5 - r/10.
\]
 There is no such $r$ for which this holds.

 The approach based on cones requires that
\[
   \frac{1}{r}  \|A_{12}\| = \frac{3}{4 r} < \|A_{22}^{-1}\|^{-1} - r \|A_{21}\|  \approx 1 - \frac{r}{10}.
\]
 We obtain
\[
   \sup r \approx 5 + \sqrt{\frac{35}{2}}.
\]
Hence we get $|\lambda| \leq 0.0817$.
\end{example}

The following example illustrates the case of the matrix $A$ for which both methods discussed above can be applied.

\begin{example}
\label{ex-gersz-our}
We put
\[
A =
\begin{bmatrix}
A_{11} & A_{12} \\
A_{21} & A_{22}
\end{bmatrix}
=\left[
\begin{array}{c|cc}
\phantom{-} 0 & \frac{2}{5} & -\frac{1}{5} \\[0.2em] \hline \\[-1em]
\phantom{-} \frac{1}{5} & \frac{3}{2} & \phantom{-} \frac{2}{5} \\[0.3em]
 -\frac{1}{10} & \frac{3}{10} & \phantom{-} 2 \\
\end{array}
\right].
\]
First, by Theorem~\ref{thm:cone-single-eigenval} we estimate the eigenvalue close to $0$. We have $a_{11}=0$ and
\begin{eqnarray*}
  \|A_{12}\|_\infty=\frac{3}{5}, \quad \|A_{21}\|_{\infty}=\frac{1}{5}, \quad \|(A_{22} - a_{11}\cdot I)^{-1}\|_\infty=\frac{5}{6}, \\
   \|(A_{22} - a_{11}\cdot I)^{-1}\|_\infty^{-2} - 4 \|A_{12}\| \cdot \|A_{21}\|_\infty= \frac{24}{25}>0.
\end{eqnarray*}
Therefore assumptions of  Theorem~\ref{thm:cone-single-eigenval} are satisfied and we obtain that the eigenvalue $\lambda$ close to $0$ satisfies $|\lambda|\leq \frac{1}{5}$.

Now we use the Gerschgorin theorems to estimate the eigenvalue close to $0$. The Gerschgorin disks are
\[
G_1(A)=\overline{B}\left(0,\frac{3}{5}\right), \quad G_2(A)=\overline{B}\left(\frac{3}{2},\frac{3}{5}\right)\quad\mbox{and}\quad
G_3(A)=\overline{B}\left(2,\frac{2}{5}\right).
\]
It is easy to see that $G_1(A)\cap G_2(A)=\emptyset$, $G_1(A)\cap G_3(A)=\emptyset$ but we rescale the matrix $A$ (with $r=3$, which is the same rescaling as in our method), we obtain the matrix
\[
\tilde{A}_r=
\begin{bmatrix}
\phantom{-} 0 & \frac{2}{15} & -\frac{1}{15} \\[0.3em]
\phantom{-} \frac{3}{5} & \frac{3}{2} & \phantom{-}\frac{2}{5} \\[0.3em]
 -\frac{3}{10} & \frac{3}{10} & \phantom{-}2 \\
\end{bmatrix},
\]
and consequently $G_1(\tilde{A}_r)\cap G_2(\tilde{A}_r)=\emptyset$ and $G_1(\tilde{A}_r)\cap G_3(\tilde{A}_r)=\emptyset$. Hence from the Gerschgorin theorem there exists an eigenvalue $\lambda$ such that $|\lambda|\leq \frac{1}{5}$.

\medskip

\textbf{Rescaling:} We look for the largest $r$ for each method, which allows us to obtain the best estimation for the eigenvalue $\lambda$ close to $0$.

For Gerschgorin disks we need to solve the following inequalities
\[
\frac{3}{2}>\frac{3}{5 r}+\frac{r}{5}+\frac{2}{5}, \quad
2>\frac{3}{5 r}+\frac{r}{10}+\frac{3}{10}.
\]
We obtain $\sup r =\frac{1}{4} \left(11+\sqrt{73}\right)$. Hence we obtain bound $|\lambda| \leq \approx 0.1228$.

The cone based approach requires
\[
 \frac{1}{r}  \|A_{12}\| = \frac{3}{5 r} < \|A_{22}^{-1}\|^{-1} - r \|A_{21}\| = \frac{6}{5} - \frac{r}{5}.
\]
 We obtain
\(
\sup r =3+\sqrt{6}.
\)
Hence we obtain the bound $|\lambda| \leq \approx 0.110102$.

By doing the same calculations as above for the transpose of the matrix $A$ we obtain
\[
\sup r = \frac{1}{2} \left(3+\sqrt{6}\right) \ \mbox{ and }\ |\lambda|\leq \approx 0.110102,
\] from the classical Gerschgorin theorem, and for cone based we get
\[
\sup r = \frac{1}{46} \left(72+\sqrt{3597}\right), \; |\lambda|\leq \approx 0.104565.
\]
As one can see the use of cone based approach gives us better estimation of the eigenvalue close to zero than the classical Gerschgorin theorem with rescaling.
\end{example}

\medskip

\textbf{Conclusions:} As one can see from above examples and theorems our method is better than Gerschgorin theorem and its modifications. The main advantages of our method are:
\begin{itemize}
\item locates the spectrum and eigenspaces of a matrix when  eigenvalues of multiplicity greater than one or clusters of very close eigenvalues are present,
\item gives better estimation for isolated eigenvalues,
\item allows to deal with composition of matrices.
\end{itemize}

\bibliographystyle{plain}
\bibliography{cone-eigenproblem}

\end{document}